\documentclass[11pt,reqno]{amsart}
\usepackage{geometry,graphicx}
\usepackage{amssymb,amsmath,amscd,graphicx,xfrac}

\usepackage{amsmath,amssymb,bbm,amscd,enumitem,amsthm}
\usepackage{graphicx,cite,multicol}

\usepackage{mathtools}
\usepackage{tikz-cd}
\usepackage{caption,subcaption}

\theoremstyle{plain}
\newtheorem{theorem}{Theorem}[section]
\newtheorem{prop}[theorem]{Proposition}

\theoremstyle{definition}

\newtheorem{definition}[theorem]{Definition}

\numberwithin{equation}{section}
 
\newcommand{\ts}{\hspace{0.5pt}}
\newcommand{\nts}{\hspace{-0.5pt}}

\DeclareMathOperator{\dens}{\mathrm{dens}}
\DeclareMathOperator{\card}{\mathrm{card}}

\DeclareMathOperator{\vol}{\mathrm{vol}}

\DeclareFontFamily{U}{mathx}{}
\DeclareFontShape{U}{mathx}{m}{n}{<-> mathx10}{}
\DeclareSymbolFont{mathx}{U}{mathx}{m}{n}
\DeclareMathAccent{\widecheck}{0}{mathx}{"71}

\newcommand{\cA}{\mathcal{A}}

\newcommand{\cL}{\mathcal{L}}
\newcommand{\vL}{\varLambda}

\newcommand{\ZZ}{\mathbb{Z}\ts}
\newcommand{\RR}{\mathbb{R}\ts}
\newcommand{\CC}{\mathbb{C}\ts}
\newcommand{\NN}{\mathbb{N}}
\newcommand{\QQ}{\mathbb{Q}}
\newcommand{\TT}{\mathbb{T}}
\newcommand{\XX}{\mathbb{X}}

\newcommand{\ii}{\ts\mathrm{i}}
\newcommand{\ee}{\,\mathrm{e}}

\newcommand{\defeq}{\mathrel{\mathop:}=}

\newcommand{\oplam}{\mbox{\Large $\curlywedge$}} 
\newcommand{\dd}{\mathop{}\!\mathrm{d}}

\newcommand{\myfrac}[2]{\frac{\raisebox{-2pt}{$#1$}}
      {\raisebox{0.5pt}{$#2$}}}
\newcommand{\bs}[1]{\boldsymbol{#1}}

\newcommand{\cvg}{\mathrm{cvg}\ts}

\title[Renormalisation for inflation systems]{Renormalisation techniques for inflation systems \\[2mm] and some of their applications}

\author{Michael Baake, Franz G\"ahler, Anna Klick, Neil Ma\~nibo, Jan Maz\'{a}\v{c}}
\address{Fakult\"at f\"ur Mathematik, Universit\"at Bielefeld, \newline
\hspace*{\parindent}Postfach 100131, 33501 Bielefeld, Germany}
\email{\{mbaake,gaehler,aklick,cmanibo,jmazac\}@math.uni-bielefeld.de }

\address{Department der Mathematik, Friedrich-Alexander-Universit\"{a}t Erlangen-N\"{u}rnberg,\newline
\hspace*{\parindent}Cauerstraße 11, 91058 Erlangen, Germany}
\email{mazacj@math.fau.de}


\begin{document} 
\begin{abstract}
Exact renormalisation techniques are important and powerful, particularly for inflation-generated systems. We review recent results in this direction. We recall the necessary notions for inflation systems and show the renormalisation principle, which allows us to obtain exact values of highly erratic functions, such as window covariograms. We apply these techniques to compute the diffraction pattern of the new monotile tilings with arbitrary precision. We also recall a recent invariant for system with pure-point spectrum, the orbit separation dimension, and its relation to renormalisation. Lastly, we recall results beyond the pure-point spectrum setting and show how  renormalisation and Lyapunov exponents can be used to exclude the presence of absolutely continuous part of the spectra. 
\end{abstract}

\maketitle

\section{Introduction and preliminaries}

The theory of aperiodic order is a mathematical discipline that emerged from the discovery of 
quasicrystals. It is concerned with the rigorous analysis of spatial structures without periodicity, 
but prevalent  long-range positional and orientational order. In view of important applications in
crystallography and materials science, one often works with tilings or Delone sets. Of particular
relevance is the frequent presence of hierarchical structures, which is often related to an underlying
inflation rule, or the the compatibility with one. All the usual suspects in two and three dimensions,
such as the Penrose and Ammann--Beenker tilings (2D), or the Ammann--Kramer and Danzer tilings (3D),
have such a hierarchical structure. Such an inflation rule gives rise to exact renormalisation
schemes for many objects, and it is the goal of this survey to demonstrate the power of such a 
structure for quantities of practical relevance. In particular, we cover correlations, covariograms,
diffraction, topological invariants, and spectra, where we employ a uniform approach via the
combination of inflation tilings with the model set description of characteristic point sets. 
We keep our exposition at an informal level, and refer to the relevant literature 
for details and proofs.

Tilings are well-studied objects in the theory of aperiodic order, which have a lot of connections to physics. A tiling can be created in various ways, for instance via inflation rules. These prescribe how to replace an inflated version of a particular type of tile (a prototile) with a collection of other tiles. The notion of an inflation rule is quite intuitive and self-explanatory, which we summarise as follows. 

\begin{definition}
Let $\{A^{}_1,\, A^{}_2,\, \dots,\, A^{}_n\}$ be a set of topologically regular, closed prototiles\index{prototile} with almost no boundary in $\RR^d$ and with finite volume. An \emph{inflation rule}\index{inflation rule} with inflation multiplier\index{inflation factor} $\lambda >1$ is a collection of mappings 
\begin{equation}
\label{eq: inflation}
\lambda A^{}_{i} \, \longmapsto \, \bigcup_{j=1}^n \bigl(A^{}_j + T^{}_{ji} \, \bigr)
\end{equation}
satisfying the following conditions:
\begin{itemize}
\item $T^{}_{ji} \subset \RR^d$ are all finite, and $A^{}_i + \varnothing = \varnothing$, 
\item the sets on the right-hand side of \eqref{eq: inflation} have pairwise disjoint interiors,
\item the image of two tiles with disjoint interiors consists of tiles with disjoint interiors, 
\item $\vol(\lambda A^{}_i) = \sum_{j=1}^n \card(T^{}_{ji}) \vol(A^{}_j)$ holds for all $1\leqslant i \leqslant n$. 
\end{itemize}
\end{definition}

Note that one can \emph{always} change the shape of the prototiles (and create new ones with a usually fractal shape called \emph{fractiles}) and turn the inflation rule \eqref{eq: inflation} into an exact equation. In such a case, we speak of a \emph{stone inflation}. 

The sets of displacements $T^{}_{ij}$  will play a crucial role in what follows. Therefore, we collect them into a~set-valued matrix $T = \bigl(T^{}_{ij}\bigr)^{}_{1\leqslant i,j \leqslant n}$, called the \emph{displacement matrix}. It appears in the literature under various names. In the case of substitution tilings with integer scaling factors, it is usually called a \emph{digit set matrix} \cite{LagariasWang,Vince}.
The matrix $M = (\card(T^{}_{ij}))$ is called the \emph{inflation matrix}\index{matrix!inflation} and takes the role of the substitution matrix. We require the inflation matrix to be primitive, which ensures that all corresponding tilings look the same locally, independent of the legal seed one starts with \cite{TAO}.

Now, in every prototile, we can select a suitable reference point, called \emph{control point}, and the entire tiling $\mathcal{T}$ can then be described as a collection of translates of prototiles. Formally,
\[
\mathcal{T} = \bigcup_{i=1}^n \vL^{}_{i} + A^{}_{i},
\]
and the sets $\vL^{}_i \subset \RR^d$ of control points inherit the inflation structure from Eq.~\eqref{eq: inflation}. This induces the equations\subjclass[2010]{
} 
\begin{equation}
\label{eq:MatrExpSys}
\vL^{}_i \, = \, \bigcup_{j=1}^{n} \,\bigl(  \lambda \vL+ t^{}_{ij} \bigr) \, = \, \bigcup_{j=1}^{n} \, \bigcup_{t \in T^{}_{ij}} 
\bigl( \lambda \vL^{}_j + t \bigr) ,
\end{equation}
for $1\leqslant i \leqslant n$, where all unions are disjoint.

An inflation tiling can only have non-trivial point spectrum (Bragg peaks) only if its inflation factor is a Pisot--Vijayaraghavan number (PV-number) \cite{Solomyak}, which is the root of an irreducible polynomial with integer coefficients, all of whose other roots are strictly smaller than $1$ in modulus. 
The same condition is also necessary for an inflation tiling to be a~cut-and-project tiling (or mutually locally derivable (MLD) with one). 
Conversely, in one dimension, the \emph{Pisot substitution conjecture} \cite{APisot} states that to have pure-point spectrum (pure Bragg diffraction), it is sufficient to have a PV-number as inflation factor and an inflation matrix with an irreducible characteristic polynomial. 
For other cases (more than one dimension, reducible characteristic polynomial), there are methods to check whether the spectrum is pure point \cite{LMS}. If this is the case, the set of control points
 $\vL^{}_{i}$ can be shown to arise from a \textit{cut-and-project scheme} (CPS)
\begin{equation}\label{eq:CPS}
  \renewcommand{\arraystretch}{1}\begin{array}{ccccc@{}l}
       \;\RR & \xleftarrow{\;\; \; \pi \;\;\; }
         & \!\RR \nts\nts \times \nts\nts H\! & 
         \xrightarrow{\;\: \pi^{}_{_\mathrm{int}} \;\: } & H & \\
         \cup & & \cup & & \cup  & \hspace*{-1ex} 
	 \raisebox{1pt}{\text{\scriptsize dense}} \\
	 \pi (\cL) & \xleftarrow{\;\ts 1-1 \;\ts } & \cL & 
	 \xrightarrow{ \qquad } &\hspace*{-2ex}\pi^{}_{_{\mathrm{int}}} (\cL) & \\
	 \| & & & & \| & \\
	 L  & \multicolumn{3}{c}{\xrightarrow{\qquad\quad\quad \  \star
		\quad\qquad\qquad}} 
         &  \,{L_{}}^{\star\nts}   &  \end{array},
  \renewcommand{\arraystretch}{1}
\end{equation}
where $\RR$ is the \emph{physical space}, a locally compact Abelian group (LCAG) $H$ is the \emph{internal space}, and $\cL$ stands for a lattice (a~co-compact, discrete subgroup) in $\RR \times H$. We also have a~pair of natural projections $\pi:\, \RR\times H \longrightarrow \RR$, $\pi^{}_{_{\mathrm{int}}}: \, \RR\times H \longrightarrow H$ such that $\pi\bigl|^{}_{\cL}$ is injective and $\pi^{}_{_{\mathrm{int}}}(\cL)$ is dense in $H$. If $H=\RR^m$, we call the CPS \emph{Euclidean}. 
The injectivity condition ensures that $\pi\bigl|^{}_{\cL}:\,\cL \longrightarrow L$ is a bijection. Therefore, one can define the \emph{star map}\index{star map} $\star : L \longrightarrow H$ for any $x\in L$ as 
\[ x \mapsto x^{\star} \defeq \pi^{}_{_{\mathrm{int}}}\bigl( (\pi|^{}_{\cL})^{-1}(x) \bigr).\]
For the CPS from above, we define a \textit{model set} as 
\[
\oplam\ts(\Omega)\,=\, \{x \in L \, : \, x^{\star}_{} \in \Omega\}\ts,
\]
where  a non-empty compact subset $\Omega \subset H$ is known as its \textit{window}. Whenever the  Pisot substitution conjecture holds, we are guaranteed that there exist windows $\Omega^{}_{i}$ such that $\vL^{}_{i} = \oplam(\Omega^{}_{i})$, up to a set of zero density, and the windows are solutions to an iterated function system. Indeed, set $\Omega^{}_{i}:= \overline{\vL^{\star}_{i}}$. Then, taking the $\star$-image and closure of Eq.~\eqref{eq:MatrExpSys} yields the desired iterated function system (IFS)
\begin{equation}\label{eq:IFS}
\Omega^{}_{i} \, = \,  \, \bigcup_{j=1}^{n} \, \bigcup_{t \in T^{}_{ij}} 
\bigl( \lambda^{\star}_{} \Omega^{}_{j} + t^{\star}_{} \bigr) \ts.
\end{equation}

The same type of embedding can be obtained whenever the inflation factor is a PV-number,
also when a pure-point spectrum is not guaranteed a priori, for instance, for tilings in $d>1$ dimensions.
We can still find a solution to \eqref{eq:IFS}, estimate its volume, and compare the densities of
the sets $\vL^{}_{i}$ and $\oplam(\varOmega^{}_{i})$. If these densities agree for all $i$,
the sets $\vL^{}_{i}$, which are always contained in $\oplam(\varOmega^{}_{i})$, are given by the
cut-and-project scheme (up to a set of density zero), so that indeed we have a system with pure-point spectrum.

For more details and further references and background, see \cite{TAO, Mazac_Thesis, SingThesis}.

\begin{figure}[!h]
    \centering
    \includegraphics[width=0.95\linewidth]{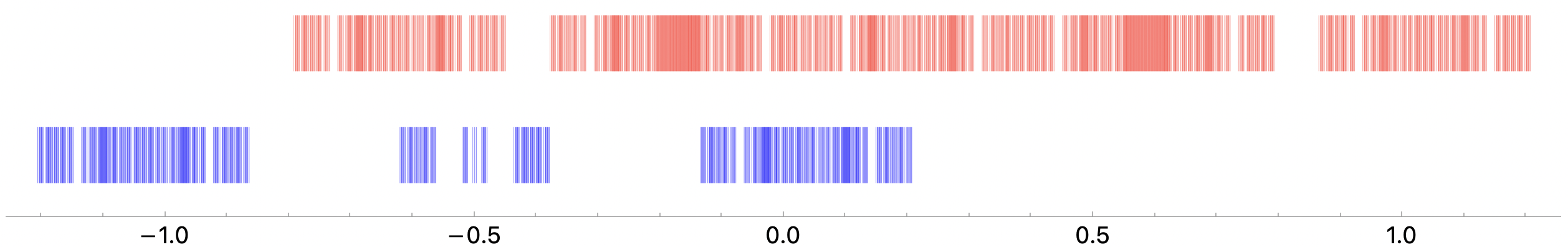}
    \caption{The window of the tiling corresponding to $\varrho$ from \eqref{eq:subst}; $\Omega^{}_{a}$ is red (top) and $\Omega^{}_{b}$ is blue (bottom). The windows are one-dimensional, but we assign some fixed arbitrary height to the points for illustration. The windows are measure-theoretically disjoint, but the resolution is limited by the large Hausdorff dimension of the window boundaries. \label{fig:ex_window} }
\end{figure}

\section{Covariograms}

\begin{figure}[!h]
    \centering
    \includegraphics[width=0.7\linewidth]{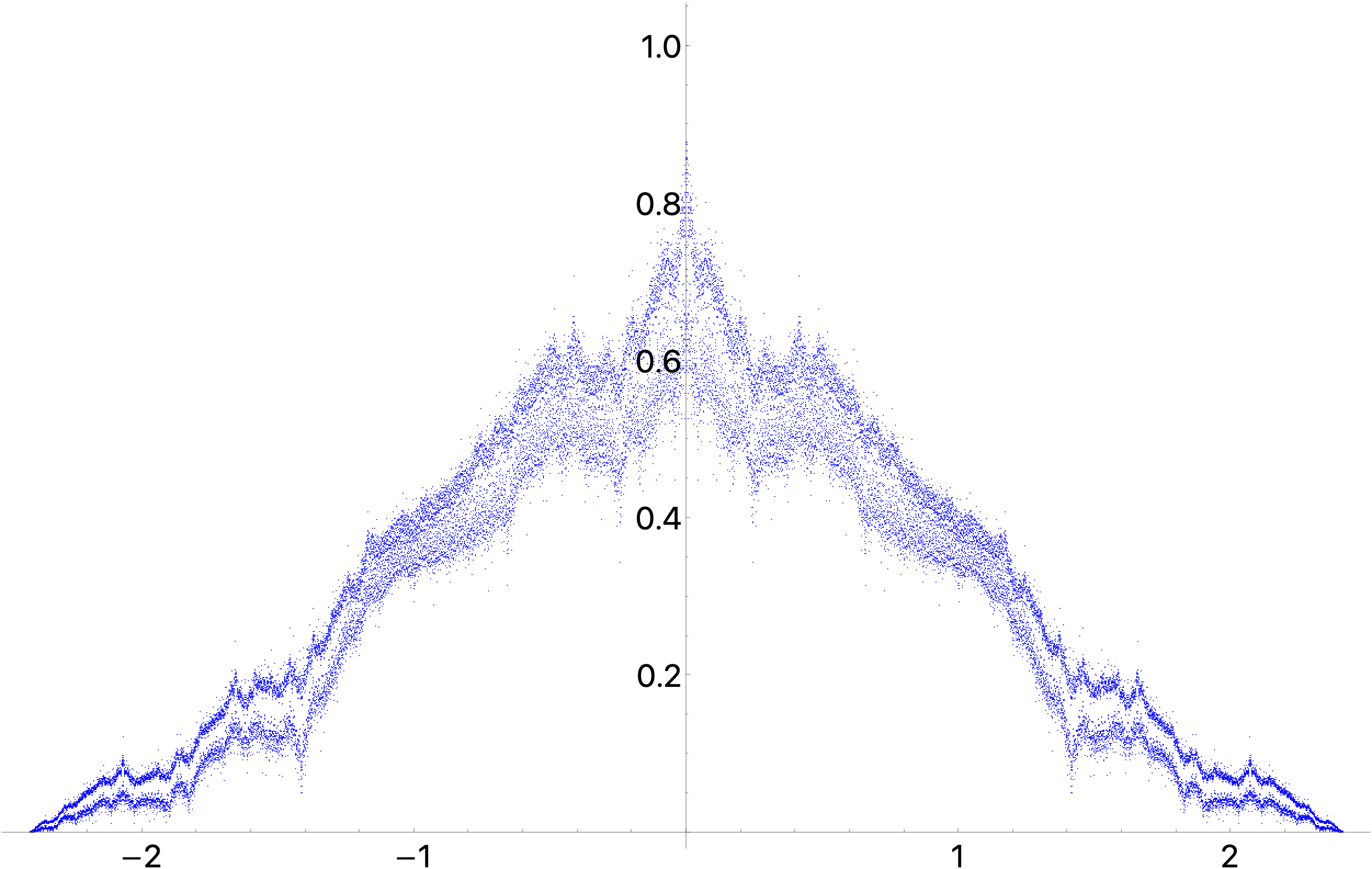}\\[3mm]\includegraphics[width=0.8\linewidth]{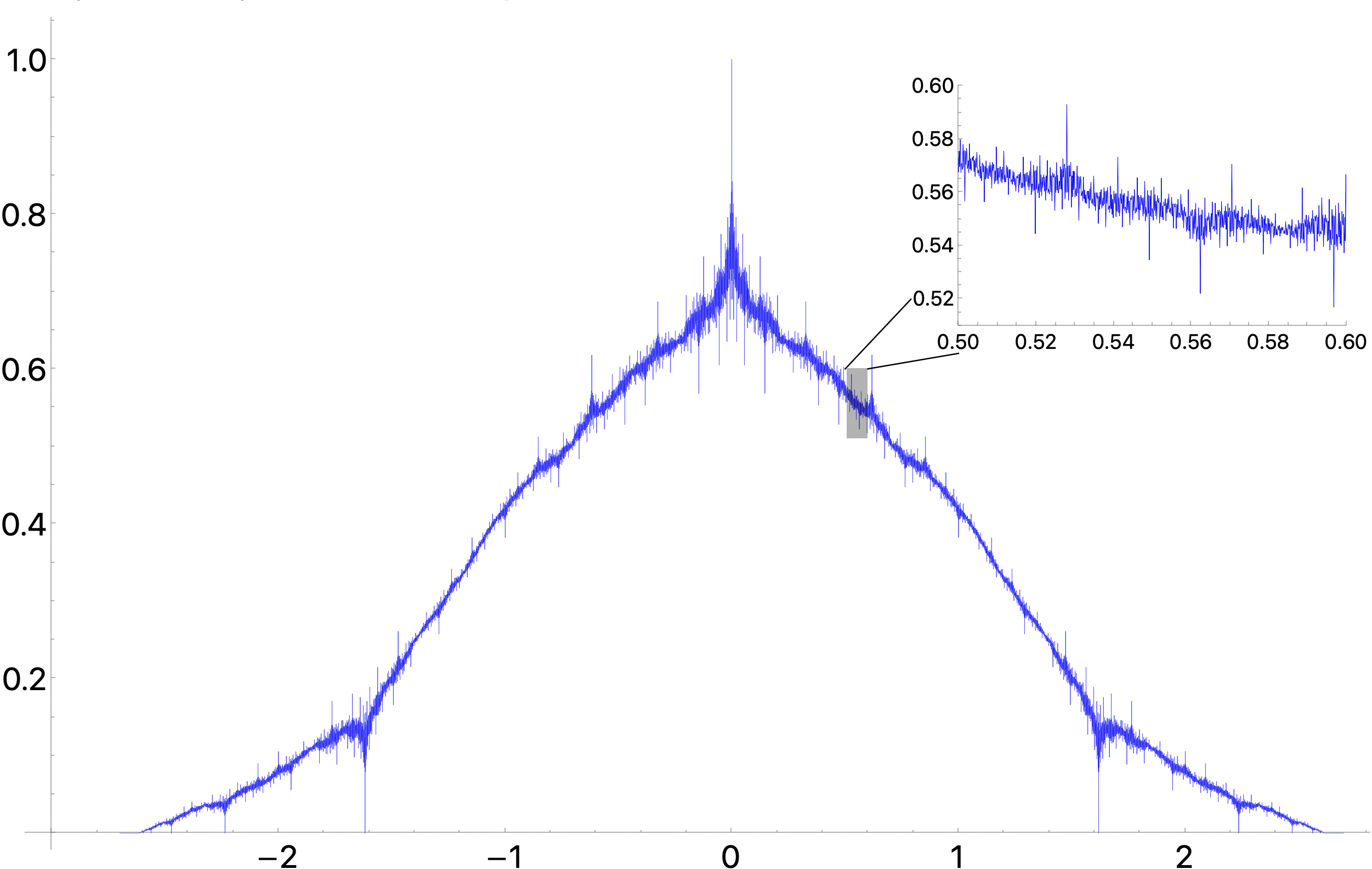}
    \caption{Upper: point plot, with $33{\ts\ts}877$ points, of the covariogram of the window of Figure~\ref{fig:ex_window}. Lower: plot with $43{\ts\ts}205$ points of the window corresponding to the reshuffled Fibonacci substitution; see \cite{BKM} for details. A small inset is included to demonstrate the highly irregular behaviour. \label{fig:cov_rfib_}}
\end{figure}

\begin{figure}[!h]
    \centering
    \includegraphics[width=0.7\linewidth]{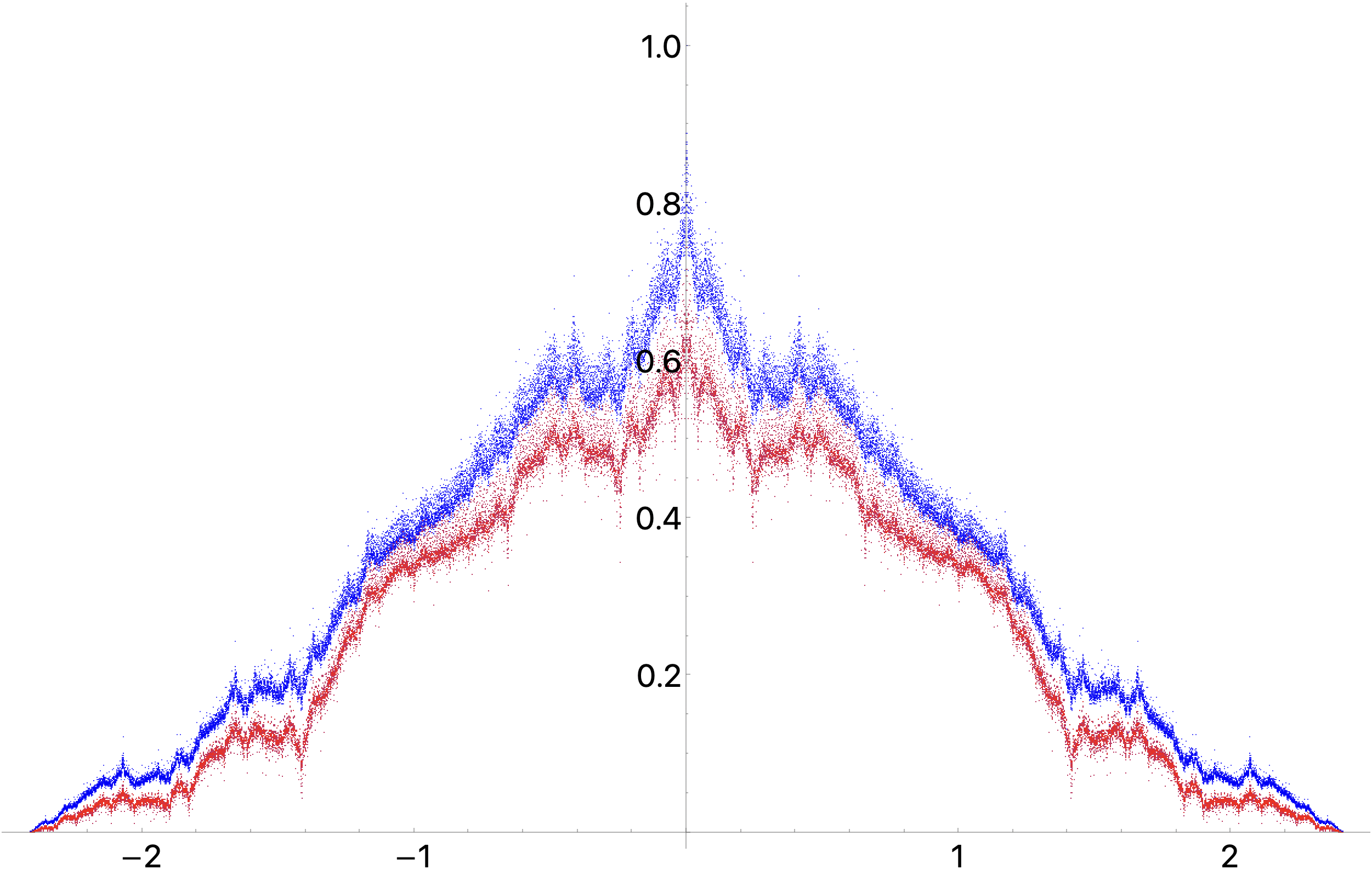}
    \caption{Point plot, with $33{\ts\ts}877$ points, of the covariogram of the window of Figure~\ref{fig:ex_window}. Here, the splitting behaviour is highlighted:  the distances $z=a\ts\lambda+b$ involving an even (odd) number of $b$'s are blue (red). \label{fig:cov_split}}
\end{figure}

The Pisot substitution conjecture holds for unimodular, primitive Pisot substitutions with irreducible characteristic polynomial on a two-letter alphabet \cite{HS03}. This permits a model set description with one-dimensional internal space $H= \RR$ in the case of unimodular substitutions. 
In general (unless the windows are finite unions of intervals), the solutions of \eqref{eq:IFS} are known as \textit{Rauzy fractals} \cite{STTopo}, but in this particular one-dimensional setting, they belong to a~special class called \textit{Cantorvals}.
\begin{definition}
A \emph{Cantorval} is a compact subset of $\RR$ with uncountably many connected components, none of which are isolated, and it is equal to the closure of its interior.
An \emph{$M$-Cantorval} (or \emph{symmetric} Cantorval) is a Cantorval with the additional property that the boundary of its interior is a Cantor set.
\end{definition}

As we will always deal with symmetric Cantorvals, we drop the adjective `symmetric' and refer to them simply as Cantorvals. These sets have a complicated structure, yet they define a regular model set, as they have positive Lebesgue measure, they are closures of their interiors and thus contain no isolated points, and have 
non-intersecting interiors and `almost no boundary'; see \cite[Cor.~6.66]{SingThesis}.

In general, one may ask when such a substitution gives rise to a Cantorval window. For example, the classic Fibonacci substitution can be described as a~model set with an interval as its window. Nevertheless, the interval is still a solution to an IFS of the form \eqref{eq:IFS}; see \cite{TAO,Klick} for details. For the class of unimodular primitive Pisot substitutions, we have the following simple criterion.

\begin{theorem}[\hspace{-0.15cm} {\cite[Thm.~4.1]{BGM24}}]
\label{thm:cantorval}
Let $\varrho$ be a primitive unimodular Pisot substitution on a binary alphabet with a model set realisation with windows $(\varOmega^{}_{a}, \, \varOmega^{}_{b})$. Denote $\varOmega = \varOmega^{}_{a} \cup \varOmega^{}_{b}$. If the boundary of the window, $\partial \varOmega$, has positive Hausdorff dimension, then $ \varOmega$, $\varOmega^{}_{a}$, and $\varOmega^{}_{b}$ are Cantorvals\index{Cantorval}.  \qed
\end{theorem}

 Moreover, if we impose some mild additional restrictions on the substitutions, we obtain their classification based on the windows; see \cite[Cor.~2.4.6]{Mazac_Thesis}. We refer the reader to \cite{BGM24, Mazac_Thesis} for a more detailed discussion on Cantorvals and their role in the theory of aperiodic order.

Now that we know when we have a Cantorval or not, the next question for us is how this influences the structure of its \textit{covariogram}.

\begin{definition}\label{def:cvg}
   Let $W \subset \RR^{m}_{}$ be a non-empty compact set. The 
   \emph{covariogram} of $W$ is 
\[
   \cvg^{}_{W}(x)\, \defeq \, \vol\ts\bigl(W \cap (x+W)\bigr) \, = \, 
   \big( \mathbf{1}^{}_{-W}  \ast \mathbf{1}^{}_{W}\big)(x)\ts,
\]
where $\mathbf{1}^{}_{W}$ is the characteristic function of $W$\hspace{-.06cm}, and~$\ast$ 
stands for the usual convolution of two functions $f,\, g\in L^{1}(\RR^{m}_{})$ given by
\[
     (f\ast g)(x) \, = \int_{\RR^m} f(y)\ts g(x-y) \dd y \ts ,
\]
with the usual understanding that this is well defined for almost all $x\in\RR^{m}_{}$.
\end{definition}

\begin{table*}[t]
\renewcommand{\arraystretch}{1.5}
\begin{center}
\begin{tabular}{ |c|c|c|c|c|c|c|c|c| } 
\hline
$z$& $0$ & $1$ & $2$ & $\lambda$ & $\lambda+1$ & $\lambda+2$ &  $2\,\lambda$ &$2\lambda+1$ \\ \hline
$\nu^{}_{aa}(z)$ & $\frac{\lambda-1}{2}$ & $0$ & $0$& $\frac{1}{2}$ & 
$\frac{5-2\lambda}{2}$ & $\frac{5-2\,\lambda}{2}$ & $\frac{9\lambda-21}{2}$ & $5-2\lambda$ \\
$\nu^{}_{ab}(z)$ & $0$ & $0$ & $0$& $\frac{\lambda-2}{2}$  & $\frac{5-2\,\lambda}{2}$ & $0$ & $\frac{22-9\lambda}{2}$ &$\frac{8\,\lambda-19}{2}$\\
$\nu^{}_{ba}(z)$ & $0$ & $\frac{\lambda-2}{2}$ & $\frac{5-2\,\lambda}{2}$& $0$ & $\frac{22-9\lambda}{2}$ & $\frac{3\,\lambda-7}{2}$ & $0$ & $\frac{3\,\lambda-7}{2}$ \\
$\nu^{}_{bb}(z)$ & $\frac{3-\lambda}{2}$ & $\frac{5-2\,\lambda}{2}$ & $0$ & $0$ & $5\lambda-12$ & $0$ & $0$ & $\frac{29-12\lambda}{2}$ \\ \hline
\end{tabular}		
\end{center}

\caption{The self-consistent part of the renormalisation equations for the tiling given by $\varrho$, with the natural tile lengths given by the left PF eigenvector.  Distances which are not possible, i.e.,~ones for which all pair correlations evaluate to $0$ due to the tile geometry, are omitted. The column corresponding to $z=0$ contains the relative tile frequencies from the (statistically normalised) right PF eigenvector of~$M^{}_{\varrho}$. \label{table:example-self-consist} }
\end{table*}

Such functions play a key role in the diffraction of model sets. First of all, the covariogram describes certain patch frequencies in our system, and it determines the autocorrelation measure of the system (up to a scaling). Next, after taking the Fourier transform (and yes, there is a rigorous way of doing this, as will be explained in the next section), the covariogram provides the intensities of the Bragg peaks. 

Now, let us consider the relative frequency of a~control point of type $i$ and a control point of type~$j$ (to the right of $i$) at a distance $x \in \RR^d$, which is given by 
\begin{equation}\label{eq:def_pair_correlation}
\nu^{}_{ij}(x) \, = \, \lim_{R\rightarrow\infty} \myfrac{\card\big(\varLambda^{(R)}_{i}\cap(\varLambda^{(R)}_{j}-x)\big)}{\card\big(\varLambda^{(R)}\big)}
\, = \, \myfrac{\dens(\varLambda^{}_i\cap (\varLambda^{}_j-x))}{\dens(\varLambda)}, 
\end{equation}
where we use the notation $A^{(R)} = A\cap B^{}_{R}(0)$ for any set $A \subseteq \RR^d$ and any $R>0$. If $\vL^{}_{i} = \oplam(\Omega^{}_{i})$, Moody's uniform distribution theorem \cite{Moody_uniform} for models sets implies
\begin{equation}
\label{eq:pair_corr_def_conv}
\nu^{}_{ij}(x) \, =  \begin{cases} \myfrac{\bigl(\boldsymbol{1}^{}_{-\varOmega^{}_{i}} \ast \boldsymbol{1}^{}_{\varOmega^{}_{j}}\bigr)(x^{\star})}{\vol(\varOmega)} \, ,  & \mbox{if} \ x \in \vL^{}_{j}{-}\vL^{}_{i}, \\
0, & \mbox{otherwise}, 
\end{cases}
\end{equation}
which gives the desired connection between covariograms and relative frequencies (or pair correlations). The inflation structure of the underlying structure (and suitable ergodic properties) permits the use of renormalisation techniques for the pair correlations as follows. 

\begin{theorem}[\hspace{-0.12cm} {\cite[Prop.~2.2.1]{NeilDiss}}]
\label{thm:renorm}
Let $\bigl(\vL^{}_{1}, \, \dots, \, \vL^{}_{n}\bigr)$ be a fixed point of a primitive geometric inflation $\varrho$ with inflation factor $\lambda >1$ arising from a substitution over an $n$-letter alphabet  $\mathcal{A}_n$. Then, the pair correlations $\nu^{}_{ij}$ exist uniformly on the hull \/ $\XX(\vL)$, and satisfy the \emph{exact} renormalisation relations 
\[
    \nu^{}_{ij}(x) \, = \, \myfrac{1}{\lambda} \sum_{k, \ell \in \cA^{}_{n}} \, \sum_{r\in T^{}_{ik}} \, \sum_{s\in T^{}_{j\ell}} \,\nu^{}_{k\ell} \left(\myfrac{x+r-s}{\lambda} \right), 
\]
where $ T = (T^{}_{ij})$ is the displacement matrix of the inflation $\varrho$. 
\end{theorem}
We note that this theorem can be extended, and using the same approach, one obtains all possible patch frequencies; see \cite{Maz25} for details.

Using Theorem~\ref{thm:renorm}, we can compute and visualise the window covariogram, as we demonstrate with an example. Consider the substitution, 
\begin{equation}
\label{eq:subst}
     \varrho \,:\, \begin{cases}
             a  \mapsto  baaaaab  \ts,\\
             \ts b \ts \mapsto  aba \ts,  \end{cases}
\end{equation}
which we abbreviate as $\varrho = (baaaaab,aba)$. Using standard techniques, we can derive its self-similar version. Since the substitution is Pisot, it possesses a model set description, with the window system shown in Figure~\ref{fig:ex_window}. It has a fractal boundary whose Hausdorff dimension is $\dim_{\mathrm{H}} (\partial W) \approx 0.91$; 
see \cite{Maz25, SingThesis, STTopo} for ways how to compute this quantity. The critical point is that calculating the covariogram via its definition as window overlaps for such a fractal is impossible. Thus, we turn to the renormalisation procedure, as described in Theorem~\ref{thm:renorm}, giving the following relations.

\begin{prop}\label{thm:ex_renorm}
     The pair correlations\/ $\nu_{\alpha \beta}$ with\/ $\alpha, \beta \in \{a,b\}$  of 
     the tiling corresponding to $\varrho$ satisfy the exact renormalisation relations 
\allowdisplaybreaks     
\begin{align*}
  \lambda^{2}_{}&\nu^{}_{aa}(z)\, = \,5\nu^{}_{aa}\bigl(\tfrac{z}{\lambda^{2}}\bigr) + 4\nu^{}_{aa}\bigl(\tfrac{z-\lambda}{\lambda^{2}}\bigr) + 4\nu^{}_{aa}\bigl(\tfrac{z+\lambda}{\lambda^{2}}\bigr) 
  +3 \nu^{}_{aa}\bigl(\tfrac{z-2\lambda}{\lambda^{2}}\bigr) + 3\nu^{}_{aa}\bigl(\tfrac{z+2\lambda}{\lambda^{2}}\bigr)+ 2\nu^{}_{aa}\bigl(\tfrac{z-3\lambda}{\lambda^{2}}\bigr) \\
  & + 2\nu^{}_{aa}\bigl(\tfrac{z+3\lambda} {\lambda^{2}}\bigr) 
   + \nu^{}_{aa}\bigl(\tfrac{z-4\lambda}{\lambda^{2}}\bigr) + \nu^{}_{aa}\bigl(\tfrac{z+4\lambda}{\lambda^{2}}\bigr)
  +\nu^{}_{ab}\bigl(\tfrac{z+1}{\lambda^{2}}\bigr)
  + \nu^{}_{ab}\bigl(\tfrac{z+\lambda+1}{\lambda^{2}}\bigr) 
  +\nu^{}_{ab}\bigl(\tfrac{z+2\lambda+1}{\lambda^{2}}\bigr)\\
  & + \nu^{}_{ab}\bigl(\tfrac{z+3\lambda+1}{\lambda^{2}}\bigr) 
  + \nu^{}_{ab}\bigl(\tfrac{z+4\lambda+1}{\lambda^{2}}\bigr) 
  + \nu^{}_{ba}\bigl(\tfrac{z-1}{\lambda^{2}}\bigr) 
   + \nu^{}_{ba}\bigl(\tfrac{z-\lambda-1}{\lambda^{2}}\bigr)
  +\nu^{}_{ba}\bigl(\tfrac{z-2\lambda-1}{\lambda^{2}}\bigr)
  + \nu^{}_{ba}\bigl(\tfrac{z-3\lambda-1}{\lambda^{2}}\bigr)\\
  &+\nu^{}_{ba}\bigl(\tfrac{z-4\lambda-1}{\lambda^{2}}\bigr) 
  +\nu^{}_{ab}\bigl(\tfrac{z-\lambda}{\lambda^{2}}\bigr) 
  + \nu^{}_{ab}\bigl(\tfrac{z}{\lambda^{2}}\bigr)
  + \nu^{}_{ab}\bigl(\tfrac{z+\lambda}{\lambda^{2}}\bigr)
  + \nu^{}_{ab}\bigl(\tfrac{z+2\lambda}{\lambda^{2}}\bigr) + \nu^{}_{ab}\bigl(\tfrac{z+3\lambda}{\lambda^{2}}\bigr)
  + \nu^{}_{ba}\bigl(\tfrac{z+\lambda}{\lambda^{2}}\bigr) \\
  &+ \nu^{}_{ba}\bigl(\tfrac{z}{\lambda^{2}}\bigr)
   +\nu^{}_{ba}\bigl(\tfrac{z-\lambda}{\lambda^{2}}\bigr) 
  + \nu^{}_{ba}\bigl(\tfrac{z-2\lambda}{\lambda^{2}}\bigr) +\nu^{}_{ba}\bigl(\tfrac{z-3\lambda}{\lambda^{2}}\bigr)
   + 2\nu^{}_{bb}\bigl(\tfrac{z}{\lambda^{2}}\bigr)
  +\nu^{}_{bb}\bigl(\tfrac{z-\lambda-1}{\lambda^{2}}\bigr) 
  +\nu^{}_{bb}\bigl(\tfrac{z+\lambda+1}{\lambda^{2}}\bigr) ,\\[4mm]
  \lambda^{2}_{}&\nu^{}_{ab}(z)\, = \, \nu^{}_{aa}\bigl(\tfrac{z-5\lambda}{\lambda^{2}}\bigr)
  + \nu^{}_{aa}\bigl(\tfrac{z-4\lambda}{\lambda^{2}}\bigr)
  +  \nu^{}_{aa}\bigl(\tfrac{z-3\lambda}{\lambda^{2}}\bigr)\\
  &+ \nu^{}_{aa}\bigl(\tfrac{z-2\lambda}{\lambda^{2}}\bigr)
  + \nu^{}_{aa}\bigl(\tfrac{z-\lambda}{\lambda^{2}}\bigr)
  + \nu^{}_{aa}\bigl(\tfrac{z+1}{\lambda^{2}}\bigr) 
  + \nu^{}_{aa}\bigl(\tfrac{z+\lambda+1}{\lambda^{2}}\bigr) 
  + \nu^{}_{aa}\bigl(\tfrac{z+2\lambda+1}{\lambda^{2}}\bigr) 
  +\nu^{}_{aa}\bigl(\tfrac{z+3\lambda+1}{\lambda^{2}}\bigr)\\
  &+ \nu^{}_{aa}\bigl(\tfrac{z+4\lambda+1}{\lambda^{2}}\bigr) 
  + \nu^{}_{ab}\bigl(\tfrac{z-\lambda+1}{\lambda^{2}}\bigr) 
  +\nu^{}_{ab}\bigl(\tfrac{z+1}{\lambda^{2}}\bigr) 
  + \nu^{}_{ab}\bigl(\tfrac{z+\lambda+1}{\lambda^{2}}\bigr)
  + \nu^{}_{ab}\bigl(\tfrac{z+2\lambda+1}{\lambda^{2}}\bigr) 
  +\nu^{}_{ab}\bigl(\tfrac{z+3\lambda+1}{\lambda^{2}}\bigr) \\
  & +\nu^{}_{ba}\bigl(\tfrac{z-5\lambda-1}{\lambda^{2}}\bigr)
  +\nu^{}_{ba}\bigl(\tfrac{z}{\lambda^{2}}\bigr) 
  +\nu^{}_{ba}\bigl(\tfrac{z+\lambda+1}{\lambda^{2}}\bigr) 
   +\nu^{}_{ba}\bigl(\tfrac{z-4\lambda}{\lambda^{2}}\bigr) 
  +\nu^{}_{bb}\bigl(\tfrac{z-\lambda}{\lambda^{2}}\bigr) 
  +\nu^{}_{bb}\bigl(\tfrac{z+1}{\lambda^{2}}\bigr)
  ,\\[4mm]
  \lambda^{2}_{}&\nu^{}_{bb}(z)\, = \,2\,\nu^{}_{aa}\bigl(\tfrac{z}{\lambda^{2}}\bigr)+ 
   \nu^{}_{aa}\bigl(\tfrac{z+5\lambda+1}{\lambda^{2}}\bigr)
   + \nu^{}_{aa}\bigl(\tfrac{z-5\lambda-1}{\lambda^{2}}\bigr)
   +  \nu^{}_{ab}\bigl(\tfrac{z-\lambda}{\lambda^{2}}\bigr)
   +\nu^{}_{ab}\bigl(\tfrac{z+4\lambda+1}{\lambda^{2}}\bigr)
   + \nu^{}_{ba}\bigl(\tfrac{z-4\lambda-1}{\lambda^{2}}\bigr)\\
   &+  \nu^{}_{ba}\bigl(\tfrac{z+\lambda}{\lambda^{2}}\bigr)
   + \nu^{}_{bb}\bigl(\tfrac{z}{\lambda^{2}}\bigr),
\end{align*}
together with\/ $\nu^{}_{ba}(z) = \nu^{}_{ab}(-z)$,  $z \in \ZZ[\sqrt{2}\,]$, where $\lambda= 1{+}\sqrt{2}$, and\/ $\nu_{\alpha \beta}(z) =0$ 
for displacements $z \notin  \vL_{\beta}- \vL_{\alpha}$.  \qed
\end{prop}

The proof of this proposition is entirely analogous to other examples included in \cite{BG16, BGM19, BKM, Klick}, and thus omitted. Now, this is an infinite set of linear equations. However, via the inflation structure, all arguments with $|z| > 2\ts\lambda$ are recursively determined from the \textit{self-consistent} part of the equations, which is given in Table \ref{table:example-self-consist}; see \cite{BGM19} for details and for a proof of the uniqueness of the solution.

With this self-consistent part in hand, serving as our seed, we can use the recursive nature of the renormalisation relations to \textit{exactly} calculate arbitrary values of the covariogram at points from $\ZZ[\sqrt{2}\ts ]$, which is a dense subset of $\RR$. 
The plot of the covariogram is shown in Figure~\ref{fig:cov_rfib_}. We see that the function seems to have a highly discontinuous nature, yet indeed is a continuous function, as it is the convolution of two functions that are both $L^{1}$ and $L^{\infty}$ \cite{Klick}. 

Moreover, this plot is an accurate representation of the function, due to the $\star$-images of points from the physical space being both uniformly distributed and dense in the window \cite{TAO}, combined with the continuity of the function. 
To account for this behaviour, a first reaction may be to attribute it to the high dimension of the window boundary; after all, this would naturally make the window overlap volume quite irregular. 
This fails to be true, as the second plot in Figure~\ref{fig:cov_rfib_}, the reshuffled Fibonacci substitution $\sigma=(aab, ba)$, comes from a window with an even higher boundary dimension of approximately $0.92$. Rather, this split, as shown in Figure~\ref{fig:cov_split}, is an artifact of the combinatorial structure of the substitution (i.e., from the number of $b$ tiles at distance $z=a\lambda+b$), which can be seen on the level of the renormalisation relations, but is outside the scope of this note (but is well worth further examination). There is no contradiction here, as the figure shows only an approximation of a continuous function. Adding more points to the plot and increasing the resolution, the gap will slowly close.

\section{Diffraction of aperiodic monotiles}

About two years ago, two families of tilings of the plane were discovered, both using a single non-convex polygon as the tile, the so-called Hat tiling \cite{Hat} and the Spectre tiling \cite{Spectre}. They both provide a partial solution to the monotile problem, as they both admit only aperiodic tilings. 

The Hat tiling consists of 12 different prototiles with respect to translations. Six Hats, which differ by a rotation by 60 degrees, and six flipped versions of Hats --- the anti-Hats, which also come in six different orientations. 
The Spectre tiling does not need the reflected version of the prototile; however, the Spectre tile appears in twelve different orientations (even though the tiling has only sixfold symmetry, as we shall see later). The twelve prototiles can be divided into two classes, Spectres versus shadow Spectres, with each tile related by a 60-degree rotation, and both classes related by a 30-degree rotation. 

For the Hat tiling as well as the Spectre tiling, there is a class of prototiles that is dominant in the tiling, meaning that the relative frequency of, say, Hats is $\tau^2$-times higher than the relative frequency of anti-Hats. Here, $\tau$ stands for the golden ratio. The same phenomenon, with a different irrational factor, happens for Spectres versus shadow Spectres. 

This small observation stands behind the reason why both tilings are considered to be only a partial solution to the monotile problem --- both tilings consist of two locally indistinguishable (LI) classes. In other words, if one places a single tile on the plane, there are precisely two types of ways to finish the tiling, and these two resulting tiling classes do not differ by a translation. In the case of Hat tiling, they are mirror images of each other, and in the Spectre case, they are related by a 30-degree rotation.

The long-range order of both tilings was determined, and it turns out that both tilings are quasiperiodic in the sense of Weyl \cite{BGS25,BGMS25}. Indeed, it was also shown that both can be obtained as a~reprojection of a model set, meaning that, for each tiling, there exists a four-dimensional Euclidean cut-and-project scheme and a window such that a non-canonical projection of the points from the strip determined by the window gives the control points of the desired tiling. These observations also imply that both tilings are pure-point diffractive \cite{Diff_Hat}. 

To derive the diffraction, i.e., the positions of the Bragg peaks and their intensities, one can use the renormalisation techniques as well (in fact, it seems to be inevitable for these tilings). To do so, one has to find a particular tiling within the Hat/Spectre family that is additionally self-similar. This can be done via a standard procedure with shape changes that build on tools from algebraic topology, see \cite{AP,CS2, Sadun} for general background, and \cite{BGS25} for a detailed treatment of the Hat tiling. This procedure yields the so-called CAP tiling (for the Hat) and CASPr tiling (for the Spectre). These tilings share the same combinatorics as the clusters of Hats/Spectres, and they are both tilings with prototiles of more than one geometric shape. Despite this, they still describe the monotile tilings in the sense that their tiling spaces are topologically conjugate with respect to the translation action of $\RR^2$. 

\begin{figure}[!h]
\centering
\includegraphics[width=0.8\textwidth]{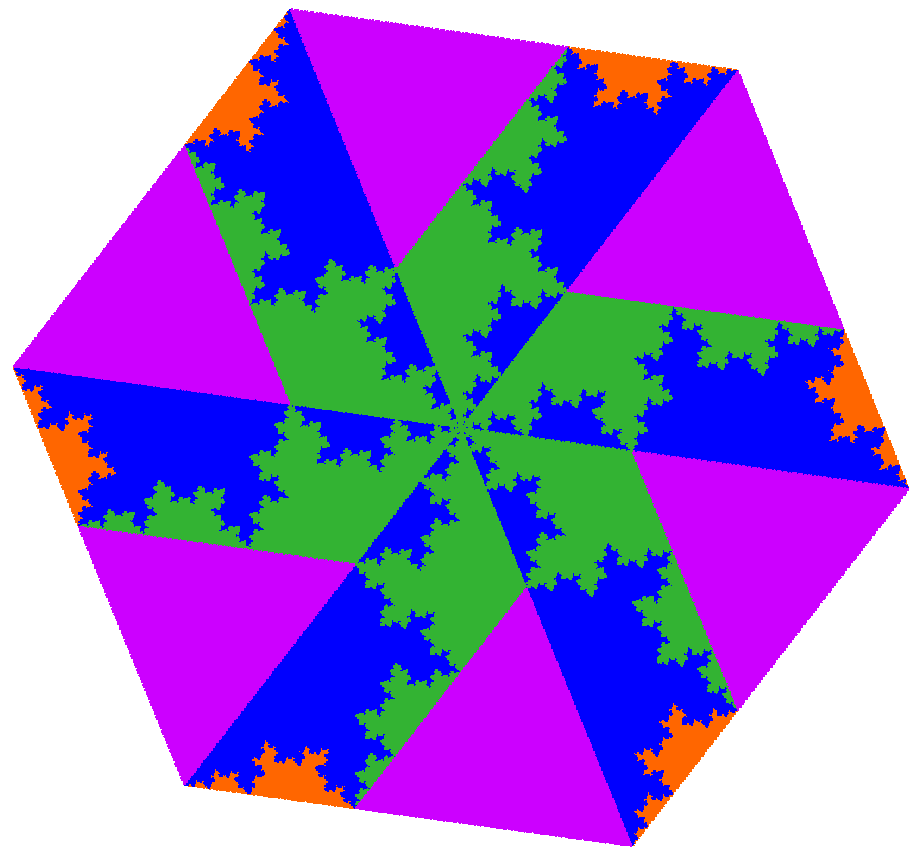}
\vspace*{-3mm}
\caption{The window for the control points of the CAP tiling. The four different colours correspond to the four different shapes of prototiles (and hence four classes of control points). For more details, see \cite{BGS25,Mazac_Thesis}.}
\label{fig:CAP-window}
\end{figure}

For the control points of both tilings, CAP and CASPr, we have equations as in \eqref{eq: inflation}. Here, the displacement matrix for the CAP tiling is of dimension 24 (4 non-congruent prototiles, each in six orientations), whereas for the CASPr, the dimension is 30 (5 non-congruent prototiles, each in six orientations).
The following theorem summarises the model set description of the control points of the CAP and the Hat tilings. For the Spectre tiling, an analogous description is possible, and we refer to \cite{Mazac_Thesis} for details.

\begin{theorem}
    The control points of the CAP tiling, $\vL^{\mathrm{CAP}}_i$, grouped according to the $24$ prototiles, are model sets with windows from Figure~$\ref{fig:CAP-window}$ that arise from a Euclidean cut-and-project scheme 
    \begin{equation*}
\renewcommand{\arraystretch}{1.2}
    \begin{array}{ccccc@{}l}
			\RR^2\simeq \CC & \hspace*{-1ex}\xleftarrow{\ \ \pi \ \ }
			& \hspace*{-2ex} \RR^2 \nts\nts \times \nts\nts \RR^2  \simeq \CC \nts\nts \times \nts\nts \CC & \hspace*{-2ex}
			\xrightarrow{\ \pi^{}_{\mathrm{int}} \ } & \hspace*{-1ex} \RR^2 \simeq \CC & \\
			\cup & & \cup & & \cup  & \hspace*{-2ex}
			\raisebox{1pt}{\scriptsize $\mathrm{dense}$} \\
			  \pi (\cL') &\hspace*{-1ex} \xleftarrow{\ \ts 1:1 \, \ } & \cL' &\hspace*{-2ex}
			\xrightarrow{ \qquad  } & \hspace*{-1ex} \pi^{}_{_{\mathrm{int}}} \nts(\cL') & \\
			\| & & & &  \| & \\ 
		      \mathcal{R}^{}_{\mathrm{CAP}} &
			\multicolumn{3}{c}{\hspace*{-1ex} \xrightarrow{\qquad\quad\quad \ \ \star \  \qquad\quad\quad}}
			&   \mathcal{R}^{\ast}_{\mathrm{CAP}}  &  
    \end{array}
\renewcommand{\arraystretch}{1}
\end{equation*}
with $\mathcal{R}^{}_{\mathrm{CAP}} = (3\tau{+}2{-}\xi)\,\ZZ[\tau,\xi]$ and the star map\index{star map} given by  
\[
(a+b\tau+c\xi+d\tau\xi)^{\star} \, = \, a+b+c+d - (b+d)\tau -(c+d)\xi + d\tau\xi,
\] with $a,b,c,d \in\QQ$, where \/$\xi$ stands for a primitive sixth root of unity. 

The set of control points of the Hat tiling, $\vL^{\mathrm{Hat}}_i$, are reprojected model sets, 
\[ \vL^{\mathrm{Hat}}_i \, = \, \bigl\{ x+ D x^{\star} \, : \, x\in \vL^{\mathrm{CAP}}_i\bigr\}\ts , \]
with 
\[ D \,= \, \myfrac{1}{16} \begin{pmatrix}
		-11 & 3\sqrt{15} \\[1mm]
		3\sqrt{15} & 11
\end{pmatrix}.\]
\end{theorem}

The diffraction theory of regular model sets is well known. For the fully Euclidean setting, the diffraction measure is pure point, and the Bragg peaks are supported by the Fourier module $L^{\circledast}$, which is the projection of the dual lattice $\cL^{\ast}$ to the physical space. For the diffraction intensities~$I$, the covariograms of the windows enter via 
\[
I(k)  = \,  \left(\myfrac{\dens(\vL)}{\vol(\varOmega)}\right)^2 \, \Big| \widehat{ \, \boldsymbol{1}^{}_{\varOmega} \ast \boldsymbol{1}^{}_{-\varOmega} \,}(-k^{\star}) \Big| = \, \left| \myfrac{\dens(\vL)}{\vol(\varOmega)} \widecheck{ \, \boldsymbol{1}^{}_{\varOmega}}(k^{\ast}) \right|^2 \, = \, |A(k)|^2, 
\]
which requires the Fourier transform of a Rauzy fractal. 
To do so, one has to use the inflation structure again and employ renormalisation techniques. For tilings with Euclidean model set description, a~method was developed by Baake and Grimm \cite{BG-Rauzy} using a matrix cocycle induced by the original inflation. The method starts with 
an $n\times n$ matrix function defined over the internal space, every entry
consisting of the inverse Fourier transform of Dirac masses placed at positions given by the $\star$-images of the entries of the displacement matrix $T$. 
For $k\in \RR^{m}$ (the internal space), the matrix elements read
\begin{equation}\label{eq:Fourier-matrix}
\underline{B}^{}_{ij}(k) \, =\, \sum_{t \in T^{}_{ij}} \ee^{2\pi\ii \langle t^{\star}  \ts | \ts k \rangle }\, ,
\end{equation} 
which is denoted by $\underline{B}(k) = \widecheck{\delta^{}_{T^{\star}}}$ for obvious reasons, and called \emph{internal Fourier matrix}.

It induces a matrix cocycle, 
\begin{equation}\label{eq:cocycle-internal}
\underline{B}^{(n)}(k) \, = \, \underline{B}(k)\,\underline{B}(\lambda^{\star}k)\, \cdots\, \underline{B}\bigl((\lambda^{\star})^{n-1}k\bigr),   
\end{equation}
which allows one to consider a matrix function $C(k)$ defined by
\begin{equation}
	\label{eq:matrix_limit}
	C(k) \, = \, \lim_{n\to\infty} \lambda^{-n}\,\underline{B}^{(n)}(k).
\end{equation}

The matrix function $C(k)$ is well defined and continuous, as the sequence $ \bigl(\lambda^{-n}\underline{B}^{(n)}(k)\bigr)^{}_{n\in\NN}$ converges compactly on $\RR^m$ \cite[Thm.~4.6]{BG-Rauzy}. Importantly, the convergence of \eqref{eq:matrix_limit} is exponentially fast, which makes it effectively computable to any desired precision. As the matrix $C(k)$ is of rank at most $1$, one can rewrite it in Dirac notation as
\[ C(k) \,= \, |c(k)\rangle \langle u |,\]
with the properly normalised left Perron--Frobenius (PF) eigenvector $\langle u|$ of the inflation matrix $M^{}_{\varrho}$. It turns out that the vector of functions $|c(k)\rangle$ has components 
\[c_i(k) \, = \, \eta\, \widecheck{\bs{1}^{}_{\,\varOmega^{}_i}}\,(k),  \] 
for some constant $\eta>0$, which can be determined explicitly, thus providing the desired quantities; see \cite[Sect.~4]{BG-Rauzy}, \cite{BG-Rauzy2} for details.

\begin{figure*}[!h]
\begin{subfigure}{.5\textwidth}
	\centering
	\includegraphics[width=0.8\textwidth]{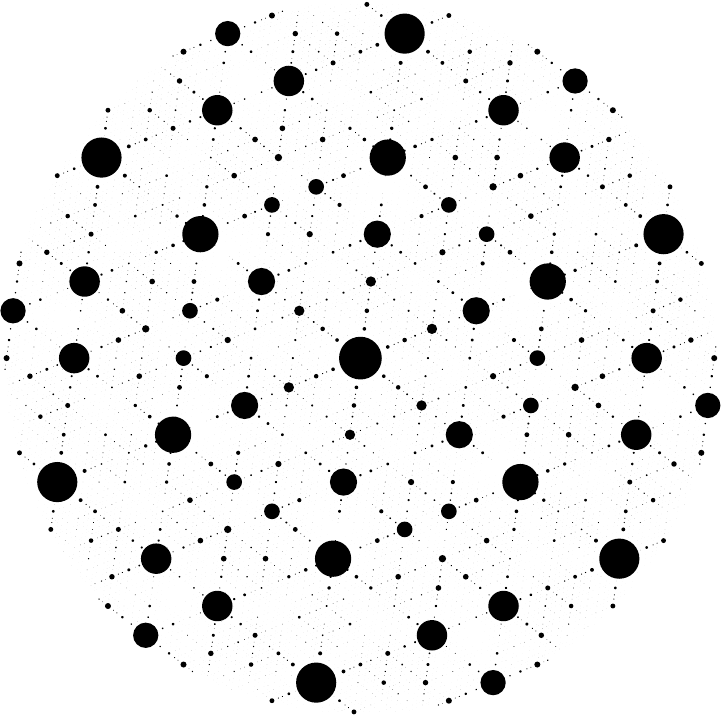}
	\caption{}
	\label{fig:CAPdiff}
	\end{subfigure}%
	\begin{subfigure}{.5\textwidth}
	\centering
	\includegraphics[width=0.8\linewidth]{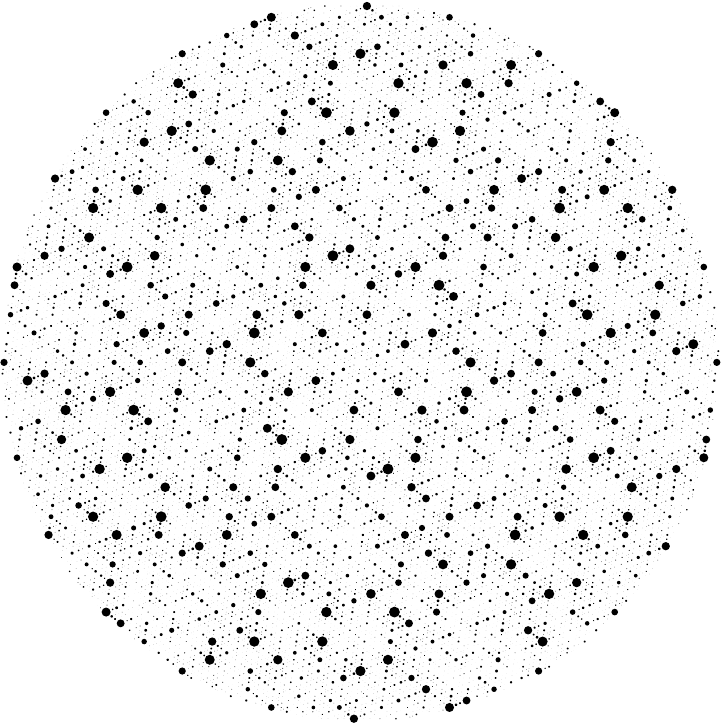}
	\caption{}
	\label{fig:CAPdiffZeroCP}
	\end{subfigure}
\vspace*{-5mm}
\caption{Diffraction pattern of the CAP tiling in the centred ball of radius~0.6. Panel (a) shows the case when all control points have equal weights, whereas Panel (b) shows the diffraction for weights chosen so that the central peak vanishes. See \cite{Diff_Hat,Mazac_Thesis} for a detailed discussion of the diffraction pattern. }
\label{fig:diffCAP}
\end{figure*}

\begin{figure*}[h!]
\begin{subfigure}{.5\textwidth}
	\centering
	\includegraphics[width=0.77\textwidth]{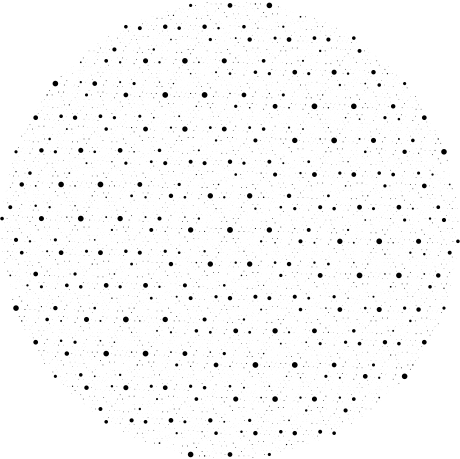}
	\end{subfigure}%
	\begin{subfigure}{.5\textwidth}
	\centering
	\scalebox{-1}[1]{\includegraphics[width=0.8\linewidth]{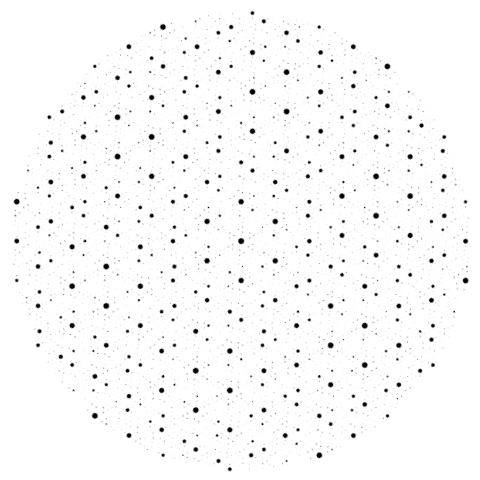}}
	\end{subfigure}
\vspace*{-5mm}
\caption{Diffraction of the Spectre tiling with equal weights for the corresponding two LI-classes of the Spectre tiling (with indicated control points), each depicted under the corresponding diffraction image (Bragg peaks in a ball of radius~0.5).  
For further details, we refer to \cite{Diff_Hat,Mazac_Thesis}. }
\label{fig:diffSpectre}
\end{figure*}

The Bragg peaks of the CAP tiling are located at the points from the Fourier module 
\[L^{\circledast}_{_{\mathrm{CAP}}} \, = \, \myfrac{(1+\xi)(\tau-\xi)\ts\ii}{3\sqrt{15}}\, \ZZ[\tau,\xi]. \]
To compute the intensities, one considers the internal Fourier matrix of dimension 24 and approximately 15 iterations of the cocycle, which is sufficient to get precision up to 10 decimal places \cite{Mazac_Thesis}. Figure \ref{fig:diffCAP} presents the diffraction image for the CAP tiling. 

\begin{figure*}[h!]
\begin{subfigure}{.5\textwidth}
	\centering
	\includegraphics[width=0.82\textwidth]{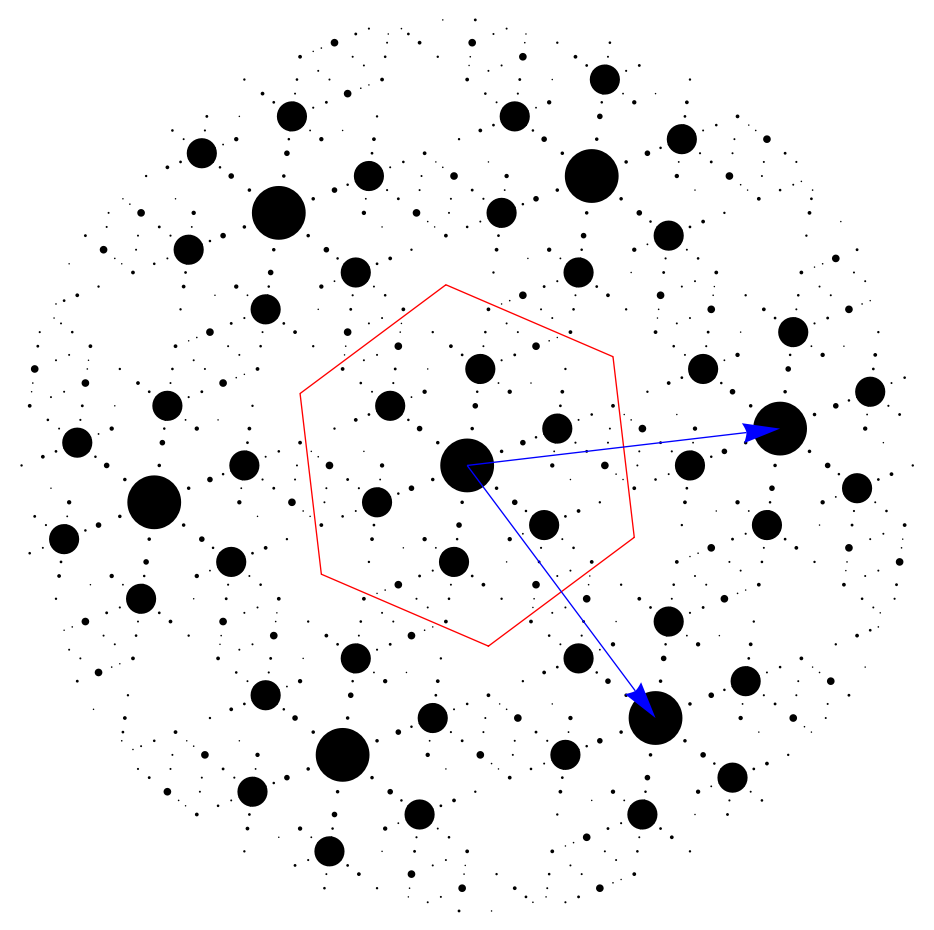}
	\end{subfigure}%
	\begin{subfigure}{.5\textwidth}
	\centering
	\scalebox{-1}[1]{\includegraphics[width=0.82\linewidth]{pics/New_Hat_Diff.png}}
	\end{subfigure}
\hfill
\begin{subfigure}{.5\textwidth}
	\centering
	\includegraphics[width=0.7\textwidth]{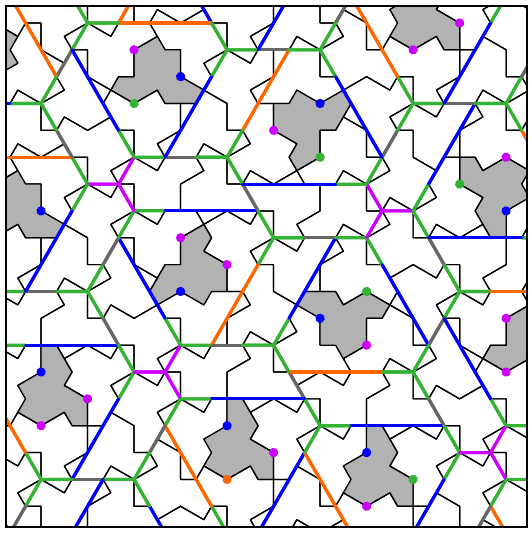}
	\end{subfigure}
	\begin{subfigure}{.5\textwidth}
	\centering
	\scalebox{-1}[1]{\includegraphics[width=0.7\linewidth]{pics/hats4jan.pdf}}
	\end{subfigure}
\vspace*{-5mm}
\caption{Diffraction of the Hat tiling with equal weights for the two LI-classes of the Hat tiling (with indicated control points), each depicted under the corresponding diffraction image. Both pictures display a lattice-periodic structure, with lattice $\tfrac{(1+\xi)(\tau-\xi)^3 \, \ii}{3\sqrt{15}}\, \ZZ[\xi]$ for the left one and its mirror image for the right one. For further details, we refer to \cite{Diff_Hat} or to \cite{Mazac_Thesis}, where a~continuous transformation from the CAP to Hat tiling on the level of the diffraction is depicted. }
\label{fig:diffHat}
\end{figure*}

To obtain the same picture for the Hat tiling, we only need to compute the diffraction intensities. Since the Hat tiling is a reprojection of the CAP tiling, this has to be reflected in the diffraction picture. Indeed, the positions of Bragg peaks remain unchanged (as the two tiling dynamical systems are topologically conjugate \cite{Diff_Hat}) and the diffraction intensities $I^{}_{\mathrm{Hat}}$ for the Hat tiling are given by the same intensity function $I^{}_{\mathrm{CAP}}$ as before, only evaluated at a modified position, namely, for any $k\in L^{\circledast}$,
\[ 
I^{}_{\mathrm{Hat}}(k)\, = \, \left| \myfrac{\dens(\vL)}{\vol( \Omega)}\,\ts \widecheck{\bs{1}^{}_{\Omega}}\bigl(k^{\star} - D^{T} (k)\bigr) \right|^2.
\]
The control points of the Hat tiling are a subset of a hexagonal lattice, which implies that, although the Hat tiling itself is aperiodic, it has a periodic diffraction pattern, with the lattice of periods being the dual lattice of the underlying hexagonal lattice (which is true for any lattice subset \cite[Thm.~10.3]{TAO}). The aperiodicity of the tiling is present in the fundamental domain, as our figures illustrate. We can use the diffraction image to distinguish the two LI classes; see Figure~\ref{fig:diffHat}.

The same applies to the Spectre tilings, whose diffraction pictures are shown in Fig.~\ref{fig:diffSpectre}. We note that, in this case, we can also use diffraction to detect the two LI-classes and to demonstrate the absence of twelvefold symmetry.

Although the approach is quite general, it will not work for the Taylor--Socolar monotile, even though it can be described as a model set. Unfortunately, one needs a non-Euclidean internal space, for which the Fourier cocycle method and the deformation arguments have not been established yet. On the other hand, as this monotile tiling is limit-periodic, a sufficiently good approximation of the diffraction image can be obtained from its periodic approximants, which is a strategy that always more or less works for limit-periodic tilings.

\section{Orbit separation dimension as topological invariant}

The \emph{orbit separation dimension} (OSD) is an invariant for dynamical systems with topological pure-point 
spectrum, which was introduced under the name \emph{amorphic complexity} \cite{FGJ16,FG20}. The intention
was to have a~tool to distinguish between different levels of
complexity for systems with zero entropy. In \cite{BGG25}, we could show that the OSD is practically computable for primitive inflation tilings, and indeed is powerful enough to distinguish many different tiling dynamical systems.

In this section, we argue that the OSD is actually closely related to autocorrelation
and diffraction. In \cite{BGG25}, it was shown that the OSD is equal to the Hausdorff
dimension of the discrete hull $\XX_0$ of the tiling (consisting of those tilings which
have a~control point at the origin), measured with a suitable pseudo-metric $D$. In this
pseudo-metric, the distance between two tilings is given by the volume fraction covered by
tiles which do not occur in both tilings, that is, the volume fraction of the region
where the two tilings do not agree. The complement is the region covered by coincident
tiles. With suitable weights for the different tile types, the latter is equal to the
total correlation between the two tilings, that is, the fraction of control points
which agree in position and type between the two tilings. Hence, the pseudo-metric
$D$ is equivalent to the autocorrelation pseudo-metric.

For inflation tilings, the metric space $(\XX_0,D)$ (or rather its projection to the \emph{maximal equicontinuous factor} (MEF), on which $D$ is a~proper metric) can be constructed as the fixed point of an
iterated function system, and its dimension is determined by the contraction rate of
the distance $D(T_1,T_2)$ of two tilings under iterated inflation. In particular,
for a tiling $T$ and its translate $T+r$ by a return vector $r$, their distance
under iterated inflation $\varrho^n$ becomes 
\[D(\varrho^n(T),\varrho^n(T)+\lambda^nr)
= D(T,T+\lambda^nr)\ts .
\]
For projection tilings, one finds
\begin{equation}\label{eq:wincorr}
  D(T,T+\lambda^nr) \, = \, \myfrac{\textrm{vol}(W \cap (W-(\lambda^nr)^\star))}{\textrm{vol}(W)},
\end{equation}
where $W$ is the window of $T$ and $(\lambda^nr)^\star$ is the $\star$-image of $\lambda^nr$.
For a Pisot inflation factor, $(\lambda^nr)^\star$ is contracting and converges to
zero, so that the two windows converge to each other. Solomyak has shown that
an inflation tiling has pure-point spectrum, and hence is a projection tiling,
if and only if $D(T,T+\lambda^nr)$ converges to zero for every return vector $r$ \cite{Solomyak}.

As we have remarked above, the OSD is determined by the contraction rate of
$ D(T,T+r)$ under inflation. This rate can be determined as follows. The
superposition of the tilings $T$ and $T+r$ is dissected into so-called overlaps,
pairs of tiles (one from each tiling) whose supports have an overlap with
non-empty interior. There are coincidence overlaps between identical tiles,
and discrepancy overlaps between non-coincident tiles. Up to translation,
there are finitely many overlap types, and the inflation on tiles induces an
inflation on overlaps. According to Solomyak, for tilings with pure-point spectrum,
every discrepancy overlap eventually produces a coincidence, so that the discrepancy
region strictly shrinks under inflation. The contraction rate is given by the
growth of the number of discrepancies (controlled by the leading eigenvalue
$\lambda^{}_{\textrm{dc}}$ of the discrepancy part of the overlap inflation),
divided by the growth of the number of all overlaps. For a tiling of $\RR^d$,
we get \cite{BGG25}
\begin{equation}\label{eq:osd}
   \textrm{OSD} \, \leqslant \,  \myfrac{\log(\lambda^d)}{\log(\lambda^d)-\log(\lambda^{}_{\textrm{dc}})},
\end{equation}
where in many cases one can actually show equality.

The same can also be obtained from \eqref{eq:wincorr}, which relates the OSD
to the rate of convergence of the sequence $(W-(\lambda^nr)^\star)$ to $W$. As one can
imagine, this depends on how complicated the boundary of $W$ is. Indeed,
also the window boundary is the fixed point of an iterated function system,
which is closely related to the overlap inflation; see \cite{Mazac_Thesis}. In nice cases (Euclidean
internal space of dimension $d^{}_\textrm{int}$, with isotropic contraction
under inflation) one finds  that
\begin{equation}\label{eq:winbd}
\begin{split}
  \textrm{OSD} \, &\leqslant \, \frac{d^{}_\textrm{int}}{d^{}_\textrm{int}-\dim(\partial W)},\\[3pt]
  \dim(\partial W) \, &\leqslant \, \myfrac{d^{}_\textrm{int}}{d}
      \myfrac{\log(\lambda^{}_{\textrm{dc}})}{\log(\lambda)},
      \end{split}
\end{equation}
where, again, one can show equality in many cases.      

We see now that the same eigenvalue $\lambda_\textrm{dc}$ governs, together with
the inflation factor $\lambda$, the convergence of the correlation between $T$
and $T+\lambda^nr$ to a limiting value, the contraction of the discrepancy
region between $T$ and $T+r$ under inflation, the dimensions of the discrete
hull $(\XX_0,D)$ and the window boundary $\partial W$, and presumably
the convergence of the Fourier--Bohr coefficients with increasing internal
space component of the $k$-vector (which depends on the dimension of the
window boundary).

As an example, we mention the Hat tiling, whose internal window boundary
(see Figure~\ref{fig:CAP-window}) has (parts of) Hausdorff dimension
\[
  \dim(\partial W) \, = \, \myfrac{\log(2+\sqrt3)}{2\log(\tau)} \, \approx \, 1.36838,
\]
whith $\tau$ as above. This relates to the OSD as
\[
  \textrm{OSD} \, = \, \myfrac{4\log(\tau)}{4\log(\tau) - \log(2+\sqrt3)} \, \approx \, 3.16644.
\]

\begin{figure}[!h]
\begin{center}
    \includegraphics[height=0.45\textwidth]{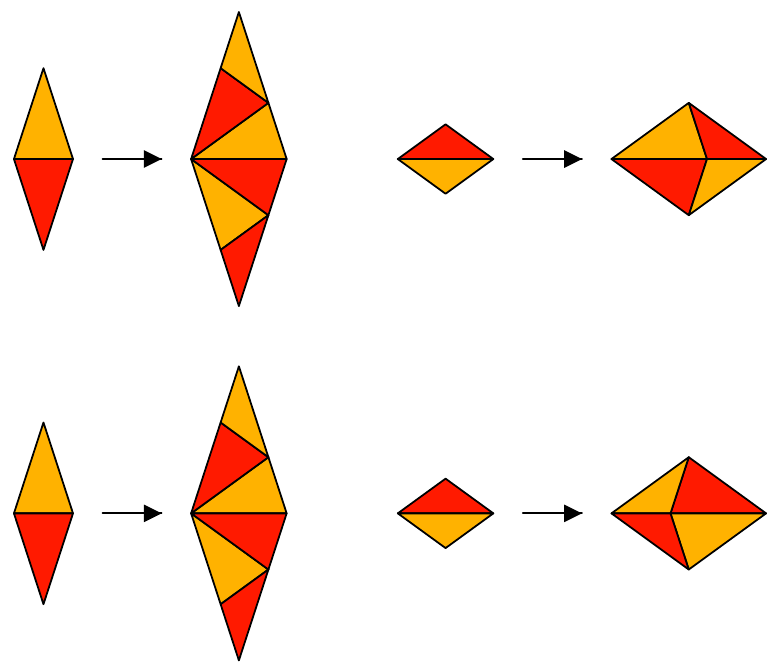}
    \end{center}
  \caption{Inflation of the Penrose tiling (top) and the pentagonal tiling
    (bottom).}
  \label{fig:penpenta}
\end{figure}

\begin{figure}[!h]
\begin{center}
    \includegraphics[height=0.45\textwidth]{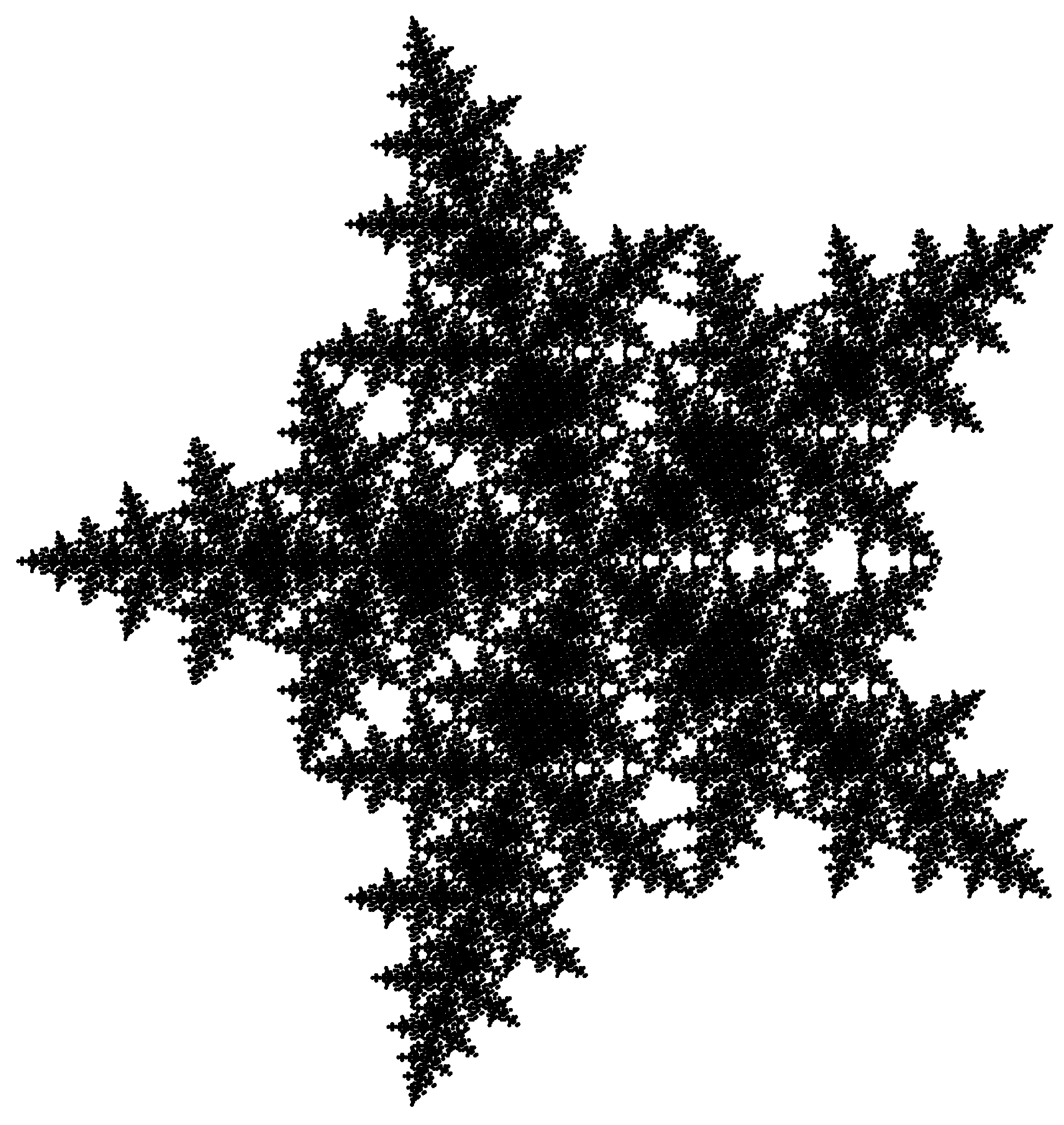}
    \end{center}
  \caption{The window of the pentagonal tiling, defined by the rule from Fig.~\ref{fig:penpenta} (bottom).}
  \label{fig:penpenta_window}
\end{figure}

As a second example, we compare a Penrose tiling (by Robinson triangles) with a
closely related pentagonal tiling (compare Figure~\ref{fig:penpenta}). For the
Penrose tiling, all four tile types occur in ten orientations each, whereas for the
pentagonal tiling, there are only five orientations per tile. The Penrose tiling
has $\textrm{OSD}=2$, and it is MLD to a pattern with a~regular decagon as its window,
whereas the pentagonal tiling has
\[
   \textrm{OSD} \, = \, \myfrac{2\log(\tau)}{2\log(\tau) - \log(1+\sqrt2)} \, \approx \, 11.8744,
\]
with a window boundary dimension
\[
   \dim(\partial W) \, = \,  \myfrac{\log(1+\sqrt2)}{\log(\tau)} \, \approx \, 1.83157.
\]
Its window is the pentagonal fractal snowflake shown in Figure~\ref{fig:penpenta_window}.

\section{Beyond pure-point diffraction --- How to exclude absolutely continuous contributions to the spectrum}

The previous sections describe the renormalisation properties satisfied by objects related to the pure-point part of the spectrum (covariogram, window, intensities of the Bragg peaks). 
Apart from these, the inflation structure also induces an exact renormalisation at the level of the \emph{absolutely continuous}
component of the diffraction, which is described by a function $h(k)\in L^{1}_{\text{loc}}(\RR^m)$ called its \emph{Radon--Nikodym density}. 

One can write $h(k)=\sum_{i,j} \overline{w^{}_i}\ts h^{}_{ij}(k) w^{}_j$, where $w^{}_i$ is the (complex) weight assigned to a control point of type $i$. The vector 
\[\boldsymbol{h}(k)=\bigl(h^{}_{11}(k), h^{}_{12}(k),\ldots,h^{}_{NN}(k)\bigr)^{T}
\]
then satisfies the equation
\begin{equation}\label{eq:ac-renorm}
\boldsymbol{h}(k)\, =\, \myfrac{1}{|\det (Q)|}P(k)\boldsymbol{h}(Q^{T}k), 
\end{equation}
where $P(k)=\overline{B}(k)\otimes B(k)$, with  $B(k)$ being the Fourier matrix in \emph{physical} space (as opposed to the internal Fourier matrix given in Eq.~\eqref{eq:Fourier-matrix}); see \cite{BGM19} for details. The main difference is that the entries are trigonometric polynomials in $\langle t \ts | \ts k \rangle$ instead of $\langle t^{\star}  \ts | \ts k \rangle$. Here, the map $Q\colon \RR^m\to \RR^m$ is the geometric expansion map associated to $\varrho$. 

One can invoke a dimension reduction argument to obtain an equation involving only $B(k)$ instead of $P(k)$. This leads to  the Fourier cocycle
\[
B^{(n)}(k)\, =\, B(k)B(Q^T k )\cdots B((Q^T)^{n-1} k),
\]
which is the analogue of Eq.~\eqref{eq:cocycle-internal} in physical space. For one-dimensional inflation rules, the inflation factor $Q=\lambda$ coincides with linear scaling. In higher dimensions, the  effective inflation factor is given via $\lambda^d=|\det(Q)|$.

The maximal \emph{Lyapunov exponent} $\chi^{B}(k)$ of this cocycle is given by 
\[
\chi^{B}(k)=\limsup_{n\to\infty} \frac{1}{n}\log\|B^{(n)}(k)\|,
\]
where the choice of norm is clearly immaterial. This quantity provides access to the nature of the absolutely continuous component, which the next result sums up; see \cite{BGM19,NeilDiss} for details. 

\begin{theorem}[\hspace{-0.15cm} {\cite[Thm.~5.7]{BGM19}}]\label{thm:absence-ac}
    Let $\varrho$ be a primitive inflation in $\RR^m$ on finitely many prototiles, with inflation map $Q$. If there exists $\varepsilon>0$ such that 
    \[
    \chi^{B}(k)\, \leqslant \, \log\sqrt{|\det(Q)|}-\varepsilon,
    \]
    for a.e.~$k\in \RR^m$, the diffraction of any (weighted) point set constructed from $\varrho$ has no absolutely continuous component. 
\end{theorem}

An additional non-degeneracy condition is assumed in the preceding theorem, namely that one has $\det(B(k))\not\equiv 0$. 
In the case where $\det(B(k))\equiv 0$, it is sometimes possible to identify the subspaces where $B(k)$ is non-invertible and restrict the analysis to the complement. 

The basic principle behind the proof is simple: if $\chi^{B}(k)$ is bounded away from $\log\sqrt{|\det(Q)|}$, one can show that $h(k)$ grows exponentially along a $Q^T$-invariant subset of $\RR^m$
of positive Lebesgue measure. On the other hand, the absolutely continuous component inherits the translation-boundedness of the entire diffraction measure. The only scenario where these two behaviours are compatible is when $h(k)=0$, for a.e. $k\in\RR^m$.

Theorem~\ref{thm:absence-ac} and the objects involved are powerful in many respects. This criterion holds in \emph{any} dimension 
and is applicable to systems with mixed spectral type (that is, those with non-trivial pure-point component). Moreover, it is able to treat systems with non-Pisot $\lambda$, all of which do not admit a covering model set, and hence do not possess a (reasonable) description via the window $W$ in internal space. 

In its current form, it is also applicable to certain systems with infinite local complexity, such as tilings that are not edge-to-edge, provided there exist only finitely many prototiles up to translation. 
We now present a procedure how to implement Theorem~\ref{thm:absence-ac} to check for singularity of the diffraction for a given example $\varrho$; see \cite{BFGR, BGM19,BGrM, NeilDiss} for fully worked-out examples. 

\begin{figure}
    \centering
    \includegraphics[height=0.45\textwidth]{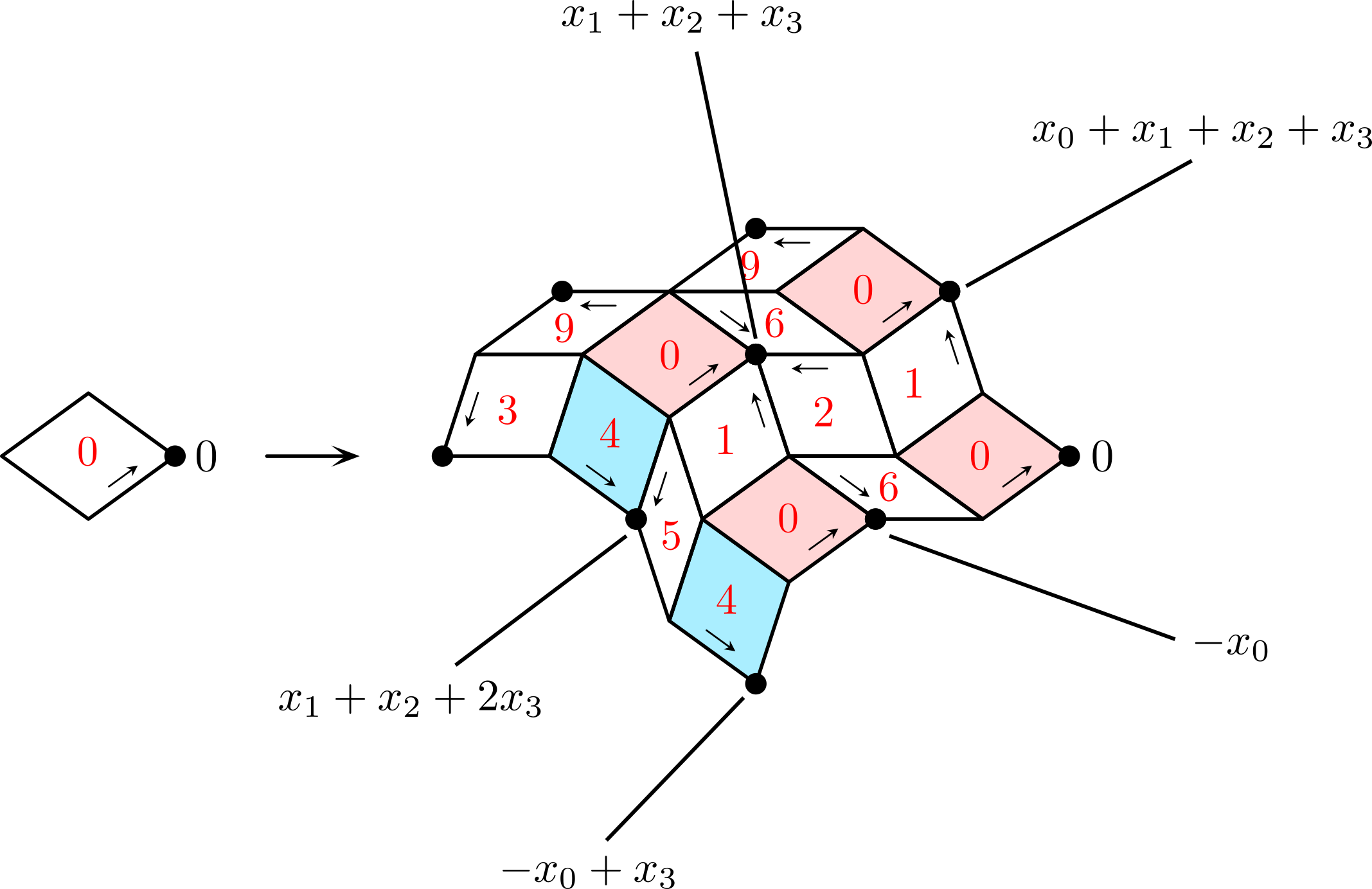}
    \caption{The level-$1$ supertile of type $0$ for the Godr\`{e}che--Lan\c{c}on--Billard (GLB) inflation rule. The location of the type-0 (pink) and type-4 (blue) tiles are labeled, where one has $x_0=\lambda\ee^{i\pi/10},\,x_{\ell}=\ee^{2\pi\ii\ell/10}x_0$ (for $1\leqslant \ell\leqslant 4$), and $\lambda=\sqrt{\frac{1}{2}(5+\sqrt{5})}$; see \cite{BGM19, NeilDiss} for a complete account.}
    \label{fig:GLB}
\end{figure}

\begin{enumerate}[leftmargin=*, label=\textbf{(\arabic*)}]
    \item From the inflation $\varrho$, compute the displacement matrix $T$.
    
    \item From the displacement matrix $T$, build the Fourier matrix $B(k)$. If the non-degeneracy condition holds for $B(k)$, one can directly build the Fourier cocycle $B^{(n)}(k)$.

    As an example, $B(k)$ for the Godr\`{e}che--Lan\c{c}on--Billard inflation rule \cite{LanBill,GL92}
    is a $10\times 10$-matrix valued function. Two of the entries are given below: 
\begin{align*}
    \big(B(k)\big)^{}_{00} &= 1+\ee^{-2\pi\ii\left\langle x_0|k\right\rangle}+\ee^{2\pi\ii\left\langle x_1+x_2+x_3|k\right\rangle} + \ee^{2\pi\ii\left\langle x_0+x_1+x_2+x_3|k\right\rangle},\\
   \big(B(k)\big)^{}_{40} &= \ee^{2\pi\ii\left\langle -x_0+x_3|k\right\rangle}+\ee^{2\pi\ii\left\langle x_1+x_2+2x_3|k\right\rangle},
\end{align*}
for $k\in \RR^2$; see Figure~\ref{fig:GLB}. 

    \item Consider $B^{(n)}(k)$ as a section of a periodic function $\widetilde{B}^{(n)}(x)$, $x\in \TT^d$ for some $d\geqslant 1$ (either via cyclotomic extensions, or via the action of the adjacency matrix of the minimal polynomial of~$\lambda$); see \cite[Sec.~5.3]{NeilDiss}.
        
    \item Compute Kingman-type bounds $K_R$ for the cocycle $\widetilde{B}^{(n)}(x)$, which are of the form
    \[
        \chi^{\widetilde{B}}(x)\, \leqslant \, \myfrac{1}{R} \int_{\TT^d}\log\|\widetilde{B}^{(R)}(x)\|\,\dd x<K_R,
    \]
    for $R\in \NN$ and for a.e. $x\in \TT^d$.
    This extends to bounds for $\chi^{B}(k)$, which is justified by sampling results for Stepanov almost periodic functions along equidistributed sequences \cite{BHL}.
        
    \item Note that Kingman-type bounds are obtained from  the integral of $\log\|\widetilde{B}^{(n)}(x)\|$ over $\TT^d$. Approximating such integrals can be done via numerical integration.
        
    \item If the obtained bound $K_R$ is less than $\log\sqrt{|\det(Q)|}$, for some $R\in \NN$, there is no absolutely continuous diffraction.
\end{enumerate}

For the GLB inflation, an upper bound $K_R$ for $\frac{1}{R}\int_{\TT^4} \log\|\widetilde{B}^{(R)}(x)\|_{\text{F}}^2\dd x$ crosses the threshold $\log\sqrt{|\det(Q)|}$ for $R=8$, where $\|\cdot\|^{ }_{\text{F}}$ is the Frobenius norm; see \cite{BGM19, NeilDiss}. 

In fact, whenever Theorem~\ref{thm:absence-ac} holds, the bounds for the Lyapunov exponent do not only provide mere absence of absolutely continuous diffraction, but also bounds for the lower local dimension of the corresponding spectral measures; see \cite{BufSol,BMS,ST}. This unlocks the route towards \emph{quantitative weak mixing}, where exponential decay rates for correlation measures (via the H\"older exponents) can be derived from the bounds for Lyapunov exponents. These are robust bounds in the sense that they are invariant under \emph{admissible deformations}, which are shape changes on the prototiles that do not alter the combinatorics of the tiling (e.g. adjacency of faces, connectedness of vertices, etc.); see \cite{ST}.

It would be interesting to apply the procedure described above to larger classes of non-Pisot systems. For example, the parametrised family of square-triangle tilings presented in \cite{SC} (originally motivated by tiling models for soft-matter quasicrystals) admits a subclass with non-Pisot inflation. Confirming singularity for this class of tilings would be a good first step in describing and visualising their spectrum.


%

\section{funding}
This work was supported by the German Research Council (Deutsche Forschungsgemeinschaft, DFG) 
under SFB-TRR 358/1 (2023) -- 491392403.




\begin{thebibliography}{99}
\itemsep=0.5pt



\bibitem{APisot}
S.~Akiyama, M.~Barge, V.~Berth{\'e}, J-Y.~Lee and A.~Siegel, 
On the Pisot substitution conjecture, in: \textit{Mathematics of Aperiodic Order}, 
eds. J.~Kellendonk, D.~Lenz and J.~Savinien,  
Birkh{\"a}user, Basel (2015), pp. 33--72.

\bibitem{AP}
J.E.~Anderson and I.F.~Putnam, Topological invariants for substitution tilings and their associated $C^*$-algebras, \textit{Ergod. Th. \& Dynam. Syst.} \textbf{18} (1998), 509--537.

\bibitem{BFGR}
M.~Baake, N.P.~Frank, U.~Grimm and E.A.~Robinson,
Geometric properties of a binary non-Pisot inflation
and absence of absolutely continuous diffraction,
\textit{Studia Math.} \textbf{247} (2019), 109--154; \texttt{arXiv:1706.03976}.

\bibitem{BG16}
M.~Baake and F.~G\"{a}hler,
Pair correlations of aperiodic inflation rules via renormalisation:
Some interesting examples,
\textit{Topol.\ Appl.} \textbf{205} (2016), 4--27; \texttt{arXiv:1511.00885}.

\bibitem{BGG25}
M.~Baake, F.~G\"ahler and P.~Gohlke, Orbit separation dimension as complexity measure for primitive inflation tilings, \textit{Ergod. Th. Dynam. Syst.} \textbf{45} (2025), 2992--3020; \texttt{arXiv:2311.03541}.

\bibitem{BGM19}
M.~Baake, F.~G\"{a}hler and N.~Ma\~{n}ibo,
Renormalisation of pair correlation measures for primitive inflation
rules and absence of absolutely continuous diffraction,
\textit{Commun.\ Math.\ Phys.} \textbf{370} (2019), 591--635; \texttt{arXiv:1805.09650}.

\bibitem{Diff_Hat}
M.~Baake, F.~G\"{a}hler, J.~Maz\'{a}\v{c} and A.~Mitchell, 
Diffraction of the Hat and Spectre tilings and some of their relatives, \textit{J. Math. Phys.} \textbf{ 66}(9) (2025), 092707:1--26; \texttt{arXiv:2502.03268}. 

\bibitem{BGMS25}
M.~Baake, F.~G\"{a}hler, J.~Maz\'{a}\v{c} and L.~Sadun,
The long-range order of the Spectre tilings, \textit{Discr. Comput. Geom.}, in press (2025);  \texttt{https://doi.org/10.1007/s00454-025-00756-z}, \texttt{arXiv:2411.15503}.
		
\bibitem{BGS25}
M.~Baake, F.~G\"{a}hler and L.~Sadun,
Dynamics and topology of the Hat family of tilings,
\textit{Israel J.\ Math.} \textbf{270} (2025), 449--485;
\texttt{arXiv:2305.05639}.

\bibitem{BGM24}
M.~Baake, A.~Gorodetski and J.~Maz\'{a}\v{c}, A naturally appearing family of Cantorvals, \textit{Lett. Math. Phys.} \textbf{114} (2024), 101:~1--11; \texttt{arXiv:2401.05372.}

\bibitem{TAO}
M.~Baake and U.~Grimm,
\textit{Aperiodic Order. Vol.\ 1: A Mathematical Invitation},
Cambridge University Press, Cambridge (2013).

\bibitem{BG-Rauzy}
M.~Baake and U.~Grimm,
Fourier transform of Rauzy fractals and point spectrum of 1D Pisot inflation tilings,
\textit{Docum.\ Math.} \textbf{25} (2020), 2303--2337;
\texttt{arXiv:1907.11012}. 

\bibitem{BG-Rauzy2}
M.~Baake and U.~Grimm,
Diffraction of a model set with complex windows, \textit{J.~Phys.: Conf. Ser.} \textbf{1458} (2020), 012006:1--7; \texttt{arXiv:1904.08285}.

\bibitem{BGrM}
M.~Baake, U.~Grimm U and N.~Ma\~nibo, 
Spectral analysis of a family of binary inflation rules,
\textit{Lett. Math. Phys.} \textbf{108} (2018), 1783--1805; \texttt{arXiv:1709.09083}.

\bibitem{BHL}
M.~Baake, A.~Haynes and D.~Lenz, Averaging almost periodic functions along exponential sequences, \newblock {In \textit{Aperiodic Order. Vol.~2:
Crystallography and Almost Periodicity}, 
M.~Baake and U.~Grimm (eds.)}, 
 \newblock Cambridge University Press, Cambridge (2017), pp.~343--362;
\texttt{arXiv:1704.08120}.
 
\bibitem{BKM} 
M.~Baake, A.~Klick and J.~Maz\'{a}\v{c}, 
Pair correlations of one-dimensional model sets and monstrous covariograms of Rauzy fractals, \textit{J. Austral. Math. Soc.} \textbf{120} (2026) 305--324; \texttt{arXiv:2502.20487}. 

\bibitem{BufSol}
A.~Bufetov and B.~Solomyak, A spectral cocycle for substitution systems and translation flows, \textit{J. Anal. Math.} \textbf{141} (2018), 165--205; \texttt{arXiv:1802.04783}.

\bibitem{BMS}
A.~Bufetov, J.~Marshall and B.~Solomyak, Local spectral estimates and quantitative weak mixing for substitution $\ZZ$-actions, \textit{J. London Math. Soc.} \textbf{111} (2025), e70136, 30~pp.

\bibitem{CS2}
A.~Clark and L.~Sadun,
When shape matters,
\textit{Ergod.\ Th.\ \& Dynam.\ Syst.} \textbf{26} (2006), 69--86;
\texttt{arXiv:math.DS/0306214}.

\bibitem{FG20}
G.~Fuhrmann and M.~Gr\"{oger},
Constant length substitutions, iterated function systems and amorphic complexity,
\textit{Math.\ Z.} \textbf{295} (2020) 1385--1404; \texttt{arXiv:1812.10789}.

\bibitem{FGJ16}
G.~Fuhrmann, M.~Gr\"{o}ger and T.~J\"{a}ger,
Amorphic Complexity,
\textit{Nonlinearity} \textbf{29} (2016) 528--565; \texttt{arXiv:1503.01036}.

\bibitem{GL92}
C.~Godr\`{e}che and F.~Lan\c{c}on, A simple example of a non-Pisot tiling with five-fold symmetry, \textit{J. Phys. I. France} \textbf{2} (1992) 207--220. 


\bibitem{HS03}
M.~Hollander and B.~Solomyak,
Two-symbol Pisot substitutions have pure discrete spectrum, \textit{Ergod.\ Th.\ \& Dynam.\ Syst.} \textbf{23} (2003), 533--540.

\bibitem{Klick}
A.~Klick, 
\textit{Averaged Shelling of Model Sets via Renormalisation}, 
Master thesis, Univ. Bielefeld (2024); available from the author.

\bibitem{LagariasWang}
J.C.~Lagarias and Y.~Wang, Integral self-affine tiles in $\RR^n$ I. Standard and nonstandard digit sets, \textit{J. Lond. Math. Soc.} \textbf{54}(1) (1996), 161--179.

\bibitem{LanBill}
F.~Lan\c{c}on and L.~Billard, Two-dimensional system with a quasi-crystalline ground state, \textit{J. Phys. France} \textbf{49} (1988), 249--256. 

\bibitem{LMS}
J.-Y.~Lee, R.V.~Moody and B.~Solomyak, 
Pure point dynamical and diffraction spectra, 
\textit{Ann. H. Poincar\'e} \textbf{3} (2002), 1003-1013; \texttt{arXiv:0910.4809}.

\bibitem{NeilDiss}
N.~Ma{\~n}ibo, \textit{Lyapunov Exponents in the Spectral Theory of Primitive Inflation Systems}, PhD Thesis, Univ. Bielefeld (2019); \texttt{urn:nbn:de:0070-pub-29359727}.

\bibitem{Mazac_Thesis}
J.~Maz\'a\v{c}, \textit{Fractal and Statistical Phenomena in Aperiodic Order}, PhD Thesis, Univ. Bielefeld (2025); \texttt{urn:nbn:de:0070-pub-30062996}. 

\bibitem{Maz25}
J.~Maz\'a\v{c}, Exact renormalisation for patch frequencies in inflation systems, \textit{preprint}, (2025); \texttt{arXiv:2507.07753}.

\bibitem{Moody_uniform}
R.V.~Moody, Uniform distribution in model sets, \textit{Can. Math. Bull.} \textbf{45} (2002), 123--130. 

\bibitem{Sadun}
L.~Sadun, \textit{Topology of Tiling Spaces}, Amer. Math. Society, Providence, RI (2008)

\bibitem{SC}
A. Say-awen and S. Coates, Octagonal tilings with three prototiles, \textit{Acta Cryst.} \textbf{A82} (2026) 179--190; \texttt{arXiv:2502.04133}.

\bibitem{STTopo}
A.~Siegel and J.~Thuswaldner, 
\textit{Topological Properties of Rauzy Fractals}, 
Soci\'{e}te Math\'{e}matique de France, Paris (2009).
 
\bibitem{SingThesis}
B.~Sing, \textit{Pisot Substitutions and Beyond}, 
PhD Thesis, Univ. Bielefeld (2007);
\texttt{urn:nbn:de:hbz:361-11555}.

\bibitem{Hat}
D.~Smith, J.S.~Myers, C.S.~Kaplan and C.~Goodman-Strauss,
An aperiodic monotile, \textit{Combin. Th.} \textbf{4}(1) (2024), 6:1--91;
\texttt{arXiv:2303.10798}.
		
\bibitem{Spectre}
D.~Smith, J.S.~Myers, C.S.~Kaplan and C.~Goodman-Strauss,
A chiral aperiodic monotile, \textit{Combin. Th.} \textbf{4}(2) (2024), 13:1--25;
\texttt{arXiv:2305.17743}.

\bibitem{Solomyak}
B.~Solomyak, Dynamics of self-similar tilings, \textit{Ergod.\ Th.\ \& Dynam.\ Syst.} \textbf{17} (1997), 695--738, and \textit{Ergod.\ Th.\ \& Dynam.\ Syst.} \textbf{19} (1999), 1685 (erratum). 

\bibitem{ST}
B.~Solomyak and R.~Trevi\~no, Spectral cocycle for substitution tilings, \textit{Ergodic Th. \& Dynam. Syst.} \textbf{44} (2024), 1629--1672; \texttt{arXiv:2201.00749}

\bibitem{Vince}
A.~Vince, Digit tiling of Euclidean space, In
\textit{Directions in Mathematical Quasicrystals},
eds.\ M.~Baake and R.~V.~Moody, Fields Institute Monographs, vol. 13,
Amer.\ Math.\ Society, Providence, RI (2000), pp. 329--370.


\end{thebibliography}

\end{document}